\documentclass[review]{elsarticle}
\usepackage{amsfonts}
\usepackage{amsmath}
\usepackage{amssymb}
\usepackage{graphicx}
\usepackage{lineno,hyperref}
\usepackage{algorithm}
\usepackage{algorithmic}
\usepackage{caption}
\usepackage{epstopdf}
\usepackage{subfigure}
\usepackage{color}
\usepackage{multirow}
\usepackage{lscape}
\usepackage{diagbox}
\usepackage{tabularx}
\usepackage{fancyhdr}
\usepackage{booktabs}

\newtheorem{assumption}{Assumption}
\newtheorem{remark}{Remark}
\newtheorem{lemma}{Lemma}
\newtheorem{theorem}{Theorem}

\newtheorem{proof}{Proof}
\modulolinenumbers[5]
\begin{document}
\captionsetup[figure]{labelfont={bf},name={Fig.},labelsep=period}
\title{(Corrected Version) The Barzilai-Borwein Method for Distributed Optimization over Unbalanced Directed Networks}
\author[a]{Jinhui Hu}
\author[b]{Xin Chen\corref{cor1}}
\ead{richard\_chen@126.com}
\author[a]{Lifeng Zheng}
\author[a]{Ling Zhang}
\author[a]{Huaqing Li\corref{cor1}}
\ead{huaqingli@swu.edu.cn}
\address[a]{Chongqing Key Laboratory of Nonlinear Circuits and Intelligent Information Processing, College of Electronic and Information Engineering, Southwest University, Chongqing $400715$, PR China}
\address[b]{School of Big Data \& Software Engineering, Chongqing University, Chongqing, $401331$ PR China.}
\cortext[cor1]{Corresponding author}

\begin{frontmatter}
\begin{abstract}
This paper studies optimization problems over multi-agent systems, in which all agents cooperatively minimize a global objective function expressed as a sum of local cost functions. Each agent in the systems uses only local computation and communication in the overall process without leaking their private information.
Based on the Barzilai-Borwein (BB) method and multi-consensus inner loops, a distributed algorithm with the availability of larger step-sizes and accelerated convergence, namely ADBB, is proposed. Moreover, owing to employing only row-stochastic weight matrices, ADBB can resolve the optimization problems over unbalanced directed networks without requiring the knowledge of neighbors’ out-degree for each agent.
Via establishing contraction relationships between the consensus error, the optimality
gap, and the gradient tracking error, ADBB is theoretically proved to converge linearly to the globally optimal solution. A real-world data set is used in simulations to validate the correctness of the theoretical analysis.
\end{abstract}
\begin{keyword}
\texttt{Multi-agent systems; distributed optimization; BB method; multi-step communications; high-performance algorithms.}
\end{keyword}
\end{frontmatter}
\thispagestyle{fancy}
\lfoot{}
\cfoot{Citation information: https://doi.org/10.1016/j.engappai.2020.104151, Engineering Applications of Artificial Intelligence. © 2020. This manuscript version is made available under the CC-BY-NC-ND 4.0 license https://creativecommons.org/licenses/by-nc-nd/4.0/.}
\renewcommand{\headrulewidth}{0mm}
\section{Introduction}
Distributed optimization is a promising paradigm that finds many practical applications such as signal processing \cite{Sun2019}, machine learning \cite{Cevher2014,Tang2020,Leottau2019}, coordinated control \cite{Xie2020,Ghasemi2020}, resource allocation \cite{Li2020h}, deep learning \cite{Lee2019,Blot2019}, and Internet of Things \cite{Liu2019}. Clearly, distributed optimization is robust than traditional centralized optimization in terms of communication networks and various practical problems can be modeled as distributed optimization, where $m$ agents cooperatively resolve the following optimization problem
\begin{equation}\label{E1}
\mathop {\min }\limits_{\tilde x \in {\mathbb{R}^n}} f\left( {\tilde x} \right) = \frac{1}{m}\sum\limits_{i = 1}^m {{f_i}\left( {\tilde x} \right)},
\end{equation}
where each agent in the system has the knowledge of only one local objective function, ${f_i}:{\mathbb{R}^n} \to \mathbb{R}$. Furthermore, each agent can only receive information from its in-neighbors and transmit information to its out-neighbors.
The mutual goal of all agents in the network is to seek the optimal solution, ${\tilde x^*}$, of problem (\ref{A1}) through communicating with its neighbors with no leak of their private information.

\subsection{Literature Review}
There is a large amount of outstanding works concerning with distributed optimization methods, including the distributed (sub)gradient method \cite{Xi2017c}, the distributed primal-dual (sub)gradient method \cite{Yuan2018a}, the distributed augmented Lagrangian method \cite{Chatzipanagiotis2015}, and the Newton method \cite{Mokhtari2017}. All distributed optimization methods can be summarized into two categories. One is Lagrangian dual variables based methods, including distributed dual decomposition \cite{Falsone2017}, distributed alternating direction method of multipliers (ADMM) \cite{Iutzeler2016a,Ling2015a,Shi2014a}. These algorithms have the superiority of achieving exact globally optimal solution while suffering from more computational complexity
than primal methods. Especially, \cite{Ling2015a,Shi2014a} are proved to converge linearly to the globally optimal solution under the conditions that the objective functions are strongly convex and have Lipschitz gradients. Some other existing methods like, DIGing \cite{Nedic2017b} and EXTRA \cite{Shi2015e}, can reach the globally optimal solution by using a sufficiently small constant step-size with no dual variables updating explicitly, which can be considered as an augmented Lagrangian primal-dual methods with a single gradient step in primal space. Another category is first-order primal methods, which are well developed in recent years owing to its simplicity and high efficiency. Some early notable methods include (sub)gradient descent (DGD) based methods \cite{Jakovetic2015,Xu2015}. In these algorithms, each agent in the system updates its local estimate according to a combination of a predefined consensus step and a local gradient descent step.

However, the common drawback of the above methods is the requirement of constructing various doubly stochastic weight matrices, which demands undirected or balanced directed networks between agents. This requirement may significantly restrict the applicability of these methods in practical applications because all agents in the network may broadcast at diverse power levels which indicates the communication capability in one direction while not in the other \cite{Xi2018a}. Therefore, some researchers are devoted to studying the algorithms over unbalanced directed networks. Based on distributed (sub)gradient descent (DGD) \cite{Jakovetic2015,Xu2015}, Xi $et$ $al.$ in \cite{Xi2017c,Xi2017b} leverage a so-called surplus-based method to realize exact convergence. Then, the (sub)gradient-push algorithm \cite{Nedic2015f} that can be applied to unbalanced directed networks constructs only column-stochastic weight matrices to achieve the globally optimal solution via incorporating a push-sum technique \cite{Xi2018d} into DGD-based methods \cite{Jakovetic2015,Xu2015}. However, the method \cite{Nedic2015f} shows relatively slow convergence rate due to the use of diminishing step-sizes. So as to accelerate the convergence, DEXTRA \cite{Xi2017h} combines the push-sum technique with EXTRA \cite{Shi2015e} to achieve linear convergence under the standard strong convexity assumption with the step-size lying in some non-trivial interval. The restriction on step-size is relaxed by the follow-up works ADD-OPT/Push-DIGing \cite{Xi2018d,Nedic2017b}, $\mathcal{A}\mathcal{B}$ \cite{Xin2018}, and $\mathcal{A}\mathcal{B}m$ \cite{Xin2019d}, where $\mathcal{A}\mathcal{B}$ simultaneously utilizes both row- and column-stochastic weight matrices to gain an accelerated convergence rate. Based on $\mathcal{A}\mathcal{B}$, $\mathcal{A}\mathcal{B}m$ first introduces distributed heavy-ball type acceleration into distributed optimization. Notice that these algorithms \cite{Xi2017b,Xi2018d,Xin2018,Xin2019,Nedic2017g} can better protect from using the doubly stochastic weight matrices and can be applied into unbalanced directed networks. However, all these methods may be impractical to some realistic environment that all agents in the networks adhere to a broadcast-based communication protocol, i.e, the agents in the network only get the knowledge of its in-degree while do not know its out-degree. Therefore, Xi $et$ $al.$ in \cite{Xi2018a} propose an elegant method which is based on a gradient tracking technique \cite{Xu2015} and uses only row-stochastic weight matrices to achieve exact convergence. FROST \cite{Xin2019} and \cite{Li2019c} extend \cite{Xi2018a} by employing uncoordinated step-sizes, which further relax the restriction on step-sizes. Theoretically, the step-sizes of FROST with linear convergence rate in the strongly convex case do not exceed $\left( {1/m{L_f}} \right)$ where $m$ is the number of agents and ${L_f}$ is the Lipschitz continuous constant. Notice that the upper bound $\left( {1/m{L_f}} \right)$ on step-sizes may be very small, which may limit the convergence rate of FROST. Thus, an useful method named Barzilai Borwein (BB) method is introduced into distributed optimization over undirected networks by DGM-BB-C \cite{Gao2022}. In \cite{Gao2022}, Gao $et$ $al.$ conduct multi-consensus inner loops to ensure both larger step-sizes and a network-independent range of step-sizes, which not only attains faster convergence than most existing works, but gives a relatively simple way of choosing the step-sizes. However, DGM-BB-C \cite{Gao2022} can only be applied into undirected or balanced directed networks due to doubly stochastic weight matrices, which may be impractical in practice.
\subsection{Motivations}
In recent literature \cite{Xi2018a,Xi2018d,Xi2017b,Xin2018,Xin2019d,Xin2019}, the notion of sufficiently small step-sizes needs to be added in theoretical results for which the agents can work in a fully distributed manner. However, researchers can always set relatively larger step-sizes in simulations, which are not in line with the theoretical results.
This paper aims to remove the notion of sufficiently small step-sizes over unbalanced directed networks and develops a distributed algorithm that converges to the globally optimal solution with fewer computational costs, communication costs and accelerated convergence.
\subsection{Contributions}
 The main contributions of this paper can be summarized as follows:

(1) A novel accelerated distributed algorithm over unbalanced directed networks constructing only row-stochastic weight matrices and using BB step-sizes, termed as ADBB, is developed to solve convex optimization problems. More significantly, ADBB is based on adapt-then-combine variation of \cite{Xi2018a} and FROST \cite{Xin2019}, and uses multi-consensus inner loops to simultaneously ensure larger step-sizes and accelerated convergence. Additionally, ADBB is theoretically demonstrated to converge to the globally optimal solution when the local objective functions are smooth and strongly convex.

(2) Compared with some recent well-known works \cite{Shi2015e,Xi2018a,Xi2018d,Xi2017h,Xin2018,Xin2019d,Nedic2017g,Xin2019}, ADBB further eliminates the restriction on step-sizes via using BB method which witnesses many successful applications in the fields of non-negative matrix factorization\cite{Huang2015a}, non-smooth optimization \cite{Huang2016}, and machine learning \cite{Tan2016b,Ma2018}. Especially, ADBB allows all agents in the network to independently automatically compute their step-sizes according to its local information, which not only relaxes the selection range of the step-sizes but also prevents from the heterogeneity problems \cite{Xin2019} existing in \cite{Nedic2017g,Xu2018b,Lu2018}.

(3) Unlike the research \cite{Xi2017c,Shi2015e,Xi2017b,Gao2022,Yuan2019b}, we consider more realistic situations, i.e., the communication networks between agents are unbalanced directed, and particularly agents exchange information in a broadcast-based directed communication network. ADBB achieves the globally optimal solution by constructing only row-stochastic weight matrices which are much easier to implement in a distributed fashion as each agent can locally decide the weights \cite{Xi2018a,Li2018}.

(4) In fact, some well-known algorithms, for example, ADD-OPT/PUSH-DIGing \cite{Xi2018d,Nedic2017b}, FROST \cite{Xin2019}, and \cite{Xi2018a,Li2019c} converge under the requirement of the largest step-sizes no exceeding $1/m{L_f}$ (possibly smaller) while ADBB increases the upper bound on the largest step-size to $1/m\mu$.
Intuitively, ADBB has more communication at each iterations. However, owing to the use of BB step-sizes and multi-consensus inner loops, ADBB converges to the globally optimal solution with the smaller number of iterations, fewer gradient-computation costs, and fewer rounds of communication than most existing works \cite{Xi2017b,Xi2018d,Xin2018,Xin2019d,Xin2019}, which is shown in Section \ref{section five}.
\subsection{Organization}
 The following organization of this paper is presented in this section. Preliminaries are presented in Section \ref{section two}. Section \ref{section three} gives the development of ADBB. The convergence of ADBB is analyzed in Section \ref{section four}. Section \ref{section five} conducts numerical experiments to resolve machine learning problems to verify the theoretical analysis. Finally, we draw a conclusion and state our future work in Section \ref{section six}.
\subsection{Basic Notations}
In this section, we present some basic notations which are useful throughout this paper.
Notice that all vectors in this paper are recognized as column vectors if no otherwise specified. The detailed definitions are given in Table \ref{Table 1}.
Based on Table \ref{Table 1}, we further introduce the Perron-vector-weighted Euclidean norm and its induced matrix norm as follows: for arbitary $\tilde x \in {\mathbb{R}^m}$ and $X \in {\mathbb{R}^{m \times m}}$
\begin{equation*}
\begin{aligned}
{\left\| {\tilde x} \right\|_\pi }: = & \sqrt {{{\left[ \pi  \right]}_1}{{\left( {{{\left[ {\tilde x} \right]}_1}} \right)}^2} + {{\left[ \pi  \right]}_2}{{\left( {{{\left[ {\tilde x} \right]}_2}} \right)}^2} +  \cdots  + {{\left[ \pi  \right]}_m}{{\left( {{{\left[ {\tilde x} \right]}_m}} \right)}^2}}  = {\left\| {{\text{diag}}\left\{ {\sqrt \pi  } \right\}\tilde x} \right\|_2},\\
{\left\| X \right\|_\pi }: = & {\left\| {{\text{diag}}\left\{ {\sqrt \pi  } \right\}X{{\left( {{\text{diag}}\left\{ {\sqrt \pi  } \right\}} \right)}^{ - 1}}} \right\|_2},
\end{aligned}
\end{equation*}
which means ${\left\| {\tilde x} \right\|_\pi } \le {{\bar \pi }^{0.5}}{\left\| {\tilde x} \right\|_2}$ and ${\left\| {\tilde x} \right\|_2} \le {{\underline{\pi } }^{ - 0.5}}{\left\| {\tilde x} \right\|_\pi }$ (see \cite{Johnson2013} for details). Denote $\vartheta  := \bar \pi /\underline{\pi }  > 1$.
\section{Preliminaries}\label{section two}
\subsection{Communication Network Model}
Consider $m$ agents exchanging information with each other over an unbalanced directed network, ${\mathcal{G}}{\text{ = }}\left( {\mathcal{V},\mathcal{E}} \right)$, where $\mathcal{V}{\text{ = }}\left\{ {1, 2, \ldots ,m} \right\}$ is the set of agents and $\mathcal{E} \subseteq \mathcal{V} \times \mathcal{V}$ is the collected ordered pairs. The weighted matrix associated with the unbalanced directed network $\mathcal{G}$ is denoted by non-negative matrix, $A = \left[ {{a_{ij}}} \right] \in {\mathbb{R}^{m \times m}}$ such that ${a_{ij}} > 0$ if $\left( {j,i} \right) \in \mathcal{E}$ and ${a_{ij}} = 0$ otherwise, where the weights ${a_{ij}}$ satisfy conditions as follows:
\begin{equation}\label{E1+}
\begin{aligned}
{a_{ij}} = \left\{ \begin{array}{l}
 > 0,\;j \in {\cal N}_i^\text{in}\\
0,\;{\rm{     otherwise}}
\end{array} \right.,\quad\sum\limits_{j = 1}^m {{a_{ij}}}  = 1,\forall i \in \mathcal{V}.
\end{aligned}
\end{equation}
Specifically, for arbitrary two agents, $i,j \in \mathcal{V}$, if $\left( {j,i} \right) \in \mathcal{E}$, then agent $j$ can transmit information to agent $i$. The in-neighbors of agent $i$ is denoted as $\mathcal{N}_i^{\text{in}}$, i.e., the set of agents can transmit messages to agent $i$. Correspondingly, the out-neighbors of agent $i$ is denoted as $\mathcal{N}_i^{\text{out}}$, i.e., the set of agents can receive information from agent $i$.
The network $\mathcal{G}$ is considered to be balanced if $\sum\nolimits_{j \in \mathcal{N}_i^{\text{out}}} {{a_{ji}} = } \sum\nolimits_{j \in \mathcal{N}_i^{\text{in}}} {{a_{ij}}} $, $i \in \mathcal{V}$, and unbalanced otherwise. Both $\mathcal{N}_i^{\text{in}}$ and $\mathcal{N}_i^{\text{out}}$ include agent $i$.
\begin{table}
\begin{tabularx}{11.5cm}{lX}  
\toprule  
\hline                    
\bf{Symbols}  & \bf{Definitions}  \\
\hline
${I_n}$   & the $n \times n$ identity matrix \\
$1_m$   & an $m$-dimensional column vector of all ones\\
${e_i}$ & an $n$-dimensional vector of all 0's except 1 at the $i$-th entry\\
${\pi}$ & the left Perron eigenvector of primitive row-stochastic matrix $A$\\
${x^{\top }}$ & transpose of vector $x$\\
${A^{\top }}$ & transpose of matrix $A$\\
${A^{ij}}$ & the $\left( {i,j} \right)$-th element of matrix $A$\\
${\left[ x \right]_i}$ & the $i$-th element of vector $x$\\
${\rm{diag}}\left\{ x \right\}$ & A diagonal matrix with all the elements of vector $x$ laying on its main diagonal\\
$X \le Y$ & each element in $Y - X$  is non-negative, where $X$ and $Y$ are two vectors or matrices\\
$X \otimes Y$ & the Kronecker product of matrices $X$ and $Y$\\
$\rho ( X )$ & the spectral radius for matrix $X$\\
$\left\|  \cdot  \right\|_2$ & the Euclidean norm for vectors and the spectral norm for matrices\\
\hline
\bottomrule
\end{tabularx}
\centering
\caption{Basic notations.}
\label{Table 1}
\end{table}
\begin{assumption}\label{A1}(\cite[Assumption 1]{Xi2018a})
The unbalanced directed network, ${\mathcal{G}}$, is strongly connected and each agent in the network has a unique identifier $i = 1, \ldots ,m$.
\end{assumption}
\begin{remark}\label{R1}
Notice that Assumption \ref{A1} is standard in recent literature \cite{Xi2018a,Xin2019,Li2019c,Li2019b}. The strongly-connected assumption on unbalanced directed networks guarantees the state averaging of the whole multi-agent system and avoids the presence of the isolated agents.
\end{remark}
\subsection{Problem Reformulation}
In this section, an equivalent form of problem (\ref{E1}) is given to help conducting the following convergence analysis.
\begin{equation}\label{E2}
\begin{aligned}
&\mathop {\min }\limits_{x \in {\mathbb{R}^{mn}}} f(x) = \frac{1}{m}\sum\limits_{i = 1}^m {{f_i}({x^i})}, \hfill \\
&{\text{s}}{\text{.t}}{\text{. }}{x^i} = {x^j}, \left( {i,j} \right) \in \mathcal{E},
\end{aligned}
\end{equation}
\begin{assumption}\label{A2}
(Smoothness): Each differentiable local objective function, $f_i$, has Lipschitz continuous gradients, i.e., for arbitrary $x,y \in {\mathbb{R}^n}$, it holds that
\begin{equation}\label{E3}
\left\| {\nabla {f_i}\left( x \right) - \nabla {f_i}\left( y \right)} \right\|_2 \le {L_f}\left\| {x - y} \right\|_2,
\end{equation}
where ${L_f} > 0$.
\end{assumption}
\begin{assumption}\label{A3}
(Strong convexity): Each local objective function, $f_i$, is strongly convex, i.e., for arbitrary $x,y \in {\mathbb{R}^n}$, it holds that
\begin{equation}\label{E4}
{f_i}\left( x \right) - {f_i}\left( y \right) \ge \left\langle {\nabla {f_i}\left( y \right),x - y} \right\rangle  + \frac{\mu }{2}{\left\| {x - y} \right\|_2^2},
\end{equation}
where $\mu > 0$.
\end{assumption}
\begin{remark}\label{R2}
Note that constants ${L_f}$, $\mu$ in Assumptions \ref{A2}-\ref{A3} satisfy $0 < \mu  \le {L_f}$ (see \cite{Bubeck2015}). Clearly, under Assumptions \ref{A3}, original problem (\ref{E1}) has a unique optimal solution which is denoted as ${{\tilde x}^*} \in {\mathbb{R}^n}$. Therefore, the globally optimal solution denoted as ${x^*} = {1_m} \otimes {{\tilde x}^*}$ to transformed problem (\ref{E2}) also exists uniquely.
We emphasize that Assumptions \ref{A2}-\ref{A3} are standard in recent literature \cite{Xi2017c,Xi2018a,Xi2018d,Xin2018,Xin2019d,Xin2019,Li2019c,Lu2018,Li2018,Li2019b}. For some specific examples, we refer the readers to \cite[Section V]{Xin2019d} and \cite[Section VI]{Xi2017b}.
\end{remark}
\section{ADBB Development}\label{section three}
\subsection{The Barzilai-Borwein Step-Sizes}
The original form of BB method solving problem (\ref{E1}) updates as follows:
\begin{equation}\label{E5}
{{\tilde x}_{k + 1}} = {{\tilde x}_k} - {\alpha _k}\nabla f\left( {{{\tilde x}_k}} \right),
\end{equation}
where ${\alpha _k}$ can be calculated by ${\alpha _k}=\left( {{\tilde s}_k^{\top }{{\tilde s}_k}} \right)/\left( {{\tilde s}_k^{\top }{{\tilde z}_k}} \right)$ or ${\alpha _k}=\left( {{\tilde s}_k^{\top }{{\tilde z}_k}} \right)/ \left( {{\tilde z}_k^{\top }{{\tilde z}_k}} \right)$ in \cite{Tan2016b}.  Therein ${{\tilde s}_k} = {{\tilde x}_k} - {{\tilde x}_{k - 1}}$ and ${{\tilde z}_k} = \nabla f\left( {{{\tilde x}_k}} \right) - \nabla f\left( {{{\tilde x}_{k - 1}}} \right)$ for $k \ge 1$. The BB method has the superiority of simplicity and flexibility. In this paper, we set the distributed BB step-sizes as follows:
\begin{equation}\label{E6}
\alpha _k^i = \frac{1}{m}\frac{{{{\left( {s_k^i} \right)}^{\top } }s_k^i}}{{{{\left( {s_k^i} \right)}^{\top } }v_k^i}},
\end{equation}
\begin{equation}\label{E7}
\alpha _k^i = \frac{1}{m}\frac{{{{\left( {s_k^i} \right)}^{\top }}v_k^i}}{{{{\left( {v_k^i} \right)}^{\top }}v_k^i}},
\end{equation}
where $s_k^i = x_k^i - x_{k - 1}^i$ and $v_k^i = \nabla {f_i}\left( {x_k^i} \right) - \nabla {f_i}\left( {x_{k - 1}^i} \right)$. One may be aware that the denominators of Eqs. \ref{E6}-\ref{E7} are tending to zero, which in fact does not affect the bounds on $\alpha _k^i$ and we give the bounds in Lemma \ref{L1}.
\begin{remark}\label{R3}
Compared with the BB step-sizes employed in DGM-BB-C \cite{Gao2022}, ADBB employs a relatively smaller step-size given in Lemma \ref{L1} ($1/m$ that of in DGM-BB-C \cite{Gao2022}) to guarantee the convergence. In fact, the division of $m$ is exactly caused by the unbalanced directed network. That is, when the unbalanced directed network degrades to undirected or balanced directed, the row-stochastic weight matrix reduces to doubly stochastic weight matrix and thus ADBB reduces to DGM-BB-C which attains larger step-sizes and a network-independent range of step-sizes.
\end{remark}
\subsection{Distributed Optimization Using Only Row-Stochastic Weight Matrices}
Distributed optimization algorithms using row-stochastic matrices are based on the method proposed in \cite{Xi2018a}, for instance, FROST \cite{Xin2019} and \cite{Li2019c} extends \cite{Xi2018a} by employing uncoordinated step-sizes. The distributed form of FROST \cite{Xin2019} updates as follows:
\begin{subequations}\label{E8}
\begin{align}
\label{E8.1}x_{k + 1}^i =& \sum\limits_{j = 1}^m {{a_{ij}}x_k^j - {\alpha _i}z_k^i},\\
\label{E8.2}y_{k + 1}^i =& \sum\limits_{j = 1}^m {{a_{ij}}y_k^j},\\
\label{E8.3}z_{k + 1}^i =& \sum\limits_{j = 1}^m {{a_{ij}}z_k^j}  + \frac{{\nabla {f_i}\left( {x_{k + 1}^i} \right)}}{{{{\left[ {y_{k + 1}^i} \right]}_i}}} - \frac{{\nabla {f_i}\left( {x_k^i} \right)}}{{{{\left[ {y_k^i} \right]}_i}}},
\end{align}
\end{subequations}
where each agent $i$ maintains three variables: $x_k^i,z_k^i \in {\mathbb{R}^{n}}$ and $y_k^i \in {\mathbb{R}^m}$, and ${{\alpha _i}}$ is the uncoordinated step-size locally selected by agent $i$, $i \in \mathcal{V}$.
Note that $a_{ij}$ in (\ref{E8}) is the weight assigned by agent $i$ to the information from agent $j$.
The row-stochastic weights, ${A^{ij}} = \left\{ {{a_{ij}}} \right\}$, are aligned with (\ref{E1+}).
\begin{remark}\label{R4}
The methods \cite{Nedic2017b,Xi2018d,Xi2017b} constructing only column-stochastic weight matrices $B = \left[ {{b_{ij}}} \right] \in {\mathbb{R}^{m \times m}}$, where
${b_{ij}} = \left\{ \begin{gathered}
   > 0,{\text{ }}i \in \mathcal{N}_j^{{\rm{out}}} \hfill \\
  0,{\text{ }}{\rm{ otherwise}} \hfill \\
\end{gathered}  \right., {\text{    }}\sum\limits_{i = 1}^m {{b_{ij}} = 1,{\text{ }}\forall j \in \mathcal{V}} $, require all agents to know their out-degrees and thus are not practical to broadcast-based communication protocol. Therefore, distributed optimization using only row-stochastic weight matrices is the state-of-the-art method, which has significant contributions than doubly or column-stochastic distributed optimization (see \cite{Xi2018a,Xin2019,Li2019b} for details). Inspired by \cite{Xi2018a,Xin2019}, this paper extends \cite{Gao2022} to unbalanced directed networks.
\end{remark}
\subsection{ADBB for Distributed Optimization}
Based on the above analysis, we develop ADBB for distributed optimization over unbalanced directed networks in Algorithm \ref{Algorithm 1}.
\begin{remark}\label{R5}
We emphasize that the initialization of auxiliary variables $y_0^i = {e_i}$, $i \in \mathcal{V}$ requires each agent in the system has a unique identifier which is satisfied under Assumption \ref{A1}. The initial values of step-sizes $\alpha _0^i$, $i \in \mathcal{V}$ should be positive and Section \ref{section five} numerically demonstrates that ADBB is not sensitive to the initial values of $\alpha _0^i$.
\end{remark}
\section{Convergence Analysis}\label{section four}
 \begin{algorithm}[!h]
\caption{ADBB for each agent $i \in \mathcal{V}$}
\label{Algorithm 1}
 \textbf{Initialization:}
For each agent $i \in \mathcal{V}$, $x_0^i \in {\mathbb{R}^n}$ is arbitrary; $y_0^i = {e_i} \in {\mathbb{R}^m}$; $z_0^i = \nabla {f_i}\left( {x_0^i} \right)\in {\mathbb{R}^n}$; $\alpha _0^i > 0$; Choose $a_{ij}$ with $\sum\nolimits_{j \in {\cal N}_i^\text{in}} {{a_{ij}}}  = 1$, for $\forall i,j \in {\cal V}$; $h = 1,2, \ldots ,H$.\\
 \textbf{For $k = 0,1, \ldots $ do} \\
 \hspace*{1em}\textbf{Each agent $i \in \mathcal{V}$ computes:} $x_{k + 1}^i\left( 0 \right) = x_k^i - \alpha _k^iz_k^i$\\
 \hspace*{1em}\textbf{Each agent $i \in \mathcal{V}$ communicates:} $x_{k + 1}^i\left( h \right) = \sum\limits_{j \in \mathcal{N}_i^{{\text{in}}} \cup \left\{ i \right\}} {{a_{ij}}x_{k + 1}^j\left( {h - 1} \right)} $\\
 \hspace*{1em} where $\alpha _k^i$ is computed by (\ref{E6}) or (\ref{E7}), and set $x_{k + 1}^i = x_{k + 1}^i\left( H \right)$.\\
 \hspace*{1em}\textbf{Each agent $i \in \mathcal{V}$ computes:} $y_{k + 1}^i\left( 0 \right) = y_k^i$\\
 \hspace*{1em}\textbf{Each agent $i \in \mathcal{V}$ communicates:} $y_{k + 1}^i\left( h \right) = \sum\limits_{j \in {\mathcal{N}_i^{{\text{in}}} } \cup \left\{ i \right\}} {{a_{ij}}y_{k + 1}^j\left( {h - 1} \right)}$\\
  \hspace*{1em} setting $y_{k + 1}^i = y_{k + 1}^i\left( H \right)$,\\
 \hspace*{1em}\textbf{Each agent $i \in \mathcal{V}$ computes:} $z_{k + 1}^i\left( 0 \right) = z_k^i + \frac{{\nabla {f_i}\left( {x_{k + 1}^i} \right)}}{{{{\left[ {y_{k + 1}^i} \right]}_i}}} - \frac{{\nabla {f_i}\left( {x_k^i} \right)}}{{{{\left[ {y_k^i} \right]}_i}}}$\\
 \hspace*{1em}\textbf{Each agent $i \in \mathcal{V}$ communicates:} $z_{k + 1}^i\left( h \right) = \sum\limits_{j \in {\mathcal{N}_i^{{\text{in}}} } \cup \left\{ i \right\}} {{a_{ij}}z_{k + 1}^j\left( {h - 1} \right)}$\\
 \hspace*{1em} and set $z_{k + 1}^i = z_{k + 1}^i\left( H \right)$.\\
 \textbf{End}
\end{algorithm}
In this section, the convergence analysis of ADBB is presented. For simplifying the analysis, we first make definitions as follows:
\begin{equation*}
\begin{aligned}
{x_k} :=& {\left[ {{{\left( {x_k^1} \right)}^{\top }}, {{\left( {x_k^2} \right)}^{\top }}, \ldots ,{{\left( {x_k^m} \right)}^{\top }}} \right]^{\top }}, \hfill \\
{{\tilde y}_k} :=& {\left[ {y_k^1, y_k^2, \ldots ,y_k^m} \right]^{\top }}, \hfill \\
{z_k} :=& {\left[ {{{\left( {z_k^1} \right)}^{\top }}, {{\left( {z_k^2} \right)}^{\top }}, \ldots ,{{\left( {z_k^m} \right)}^{\top }}} \right]^{\top }}, \hfill \\
{y_k} :=& {{\tilde y}_k} \otimes {I_n}, \hfill \\
{\hat y_k} :=& {\text{diag}}\left\{ {{y_k}} \right\}, \hfill \\
{\alpha _k} :=& {\left[ {\alpha _k^1, \alpha _k^2, \ldots ,\alpha _k^m} \right]^{\top }},  \hfill \\
D_k^\alpha  :=& {\rm{diag}}\left\{ {{\alpha _k}} \right\} \otimes {I_n}, \hfill \\
\mathcal{A} :=& A \otimes {I_n} \in {\mathbb{R}^{mn \times mn}}, \hfill \\
\nabla F\left( {{x_k}} \right) :=& {\left[ {\nabla {f_1}{{\left( {x_k^1} \right)}^{\top }}, \nabla {f_2}{{\left( {x_k^2} \right)}^{\top }}, \ldots ,\nabla {f_m}{{\left( {x_k^m} \right)}^{\top }}} \right]^{\top }}, \hfill\\
\nabla F\left( {{x^*}} \right) :=& {\left[ {\nabla {f_1}{{\left( {{x^*}} \right)}^{\top }}, \nabla {f_2}{{\left( {{x^*}} \right)}^{\top }}, \ldots ,\nabla {f_m}{{\left( {{x^*}} \right)}^{\top }}} \right]^{\top }}, \hfill
\end{aligned}
\end{equation*}
where ${x_k},{\tilde y_k},{z_k}, {\alpha _k},\nabla F\left( {{x_k}} \right),\nabla F\left( {{x^*}} \right)$, simultaneously collect their local variables.
Then, based on the above definitions, we give the compact form of ADBB as follows:
\begin{subequations}\label{E10}
\begin{align}
\label{E10.1} {x_{k + 1}} = & {{\cal A}^H}\left( {{x_k} - D_k^\alpha {z_k}} \right),\\
\label{E10.2} {\tilde y_{k + 1}} = & {A^H}{\tilde y_k},\\
\label{E10.3} {z_{k + 1}} = & {{\cal A}^H}\left( {{z_k} + \hat y_{k + 1}^{ - 1}\nabla F\left( {{x_{k + 1}}} \right) - \hat y_k^{ - 1}\nabla F\left( {{x_k}} \right)} \right),
\end{align}
\end{subequations}
where ${\tilde y_0} = {I_m}$, ${z_0} = \nabla f\left( {{x_0}} \right)$ and ${x_0}$ is arbitrary. This paper extends one dimension in \cite{Xi2018a} to $n$ dimensions.
To proceed, we continue to define
\begin{equation*}
\begin{aligned}
 {y_\infty } :=& \mathop {\lim }\limits_{k \to \infty } {y_k}, \hfill \\
{\hat y_\infty } :=& {\rm{diag}}\left\{ {{y_\infty }} \right\}, \hfill \\
Y :=& \mathop {\sup }\limits_{k \ge 0} \left\| {\hat y_k^{ - 1}} \right\|_2, \hfill \\
\hat Y :=& \mathop {\sup }\limits_{k \ge 0} \left\| {{y_k}} \right\|_2, \hfill \\
{p_1} :=& \left\| {{I_{mn}} - {\mathcal{A}^H}} \right\|_2, \hfill \\
\end{aligned}
\end{equation*}
 \begin{equation*}
 \begin{aligned}
 {\alpha _{\max }} :=& \mathop {\max }\limits_{k \ge 0} \left\{ {\alpha _k^i} \right\}, \hfill \\
 {{\bar \alpha }_{\max }} :=& \mathop {\max }\limits_{k \ge 0} \left\{ {\frac{1}{m}\sum\limits_{i = 1}^m {\alpha _k^i} } \right\}, \hfill
 \end{aligned}
 \end{equation*}
where ${y_\infty },Y,\hat Y$ exist due to the primitive row-stochastic matrix $A$.

\subsection{Auxiliary Relations}
\begin{lemma}\label{L1}
Supposing that Assumptions \ref{A2}-\ref{A3} hold, for $\forall k \ge 0$, the BB step-size $\alpha _k^i$, $i \in \mathcal{V}$ computed by (\ref{E6}) or (\ref{E7}) in ADBB satisfies
\begin{equation}\label{E11}
\frac{1}{{m{L_f}}} \le \alpha _k^i \le \frac{1}{{m\mu }}.
\end{equation}
\begin{proof}
To begin with, we derive bounds on the BB step-size ${\alpha _k^i}$ following (\ref{E6}). According to the strong convexity of local objective function ${f_i}$, it holds that
\begin{equation}
{\left( {x_k^i - x_{k - 1}^i} \right)^{\top }}\left( {\nabla{f_i}\left( {x_k^i} \right) - \nabla{f_i}\left( {x_{k - 1}^i} \right)} \right) \ge \mu {\left\| {x_k^i - x_{k - 1}^i} \right\|_2^2}.
\end{equation}
Then, the upper bound for each BB step-size ${\alpha _k^i}$ is derived as follows:
\begin{equation}
\begin{aligned}
\alpha_k^i =& \frac{1}{m}\frac{{{{\left( {x_k^i - x_{k - 1}^i} \right)}^{\top }}\left( {x_k^i - x_{k - 1}^i} \right)}}{{{{\left( {x_k^i - x_{k - 1}^i} \right)}^{\top }}\left( {\nabla {f_i}\left( {x_k^i} \right) - \nabla {f_i}\left( {x_{k - 1}^i} \right)} \right)}}\\
 \le& \frac{1}{m}\frac{{{{\left\| {x_k^i - x_{k - 1}^i} \right\|}_2^2}}}{{\mu {{\left\| {x_k^i - x_{k - 1}^i} \right\|}_2^2}}}\\
 =&\frac{1}{m \mu },
 \end{aligned}
\end{equation}
and according to Lipschitz continuity and Cauchy inequality, the lower bound for each BB step-size ${\alpha _k^i}$ is derived as follows:
\begin{equation}
\begin{aligned}
\alpha _k^i =& \frac{1}{m}\frac{{{{\left( {x_k^i - x_{k - 1}^i} \right)}^{\top }}\left( {x_k^i - x_{k - 1}^i} \right)}}{{{{\left( {x_k^i - x_{k - 1}^i} \right)}^{\top }}\left( {\nabla {f_i}\left( {x_k^i} \right) - \nabla {f_i}\left( {x_{k - 1}^i} \right)} \right)}}\\
 \ge& \frac{1}{m}\frac{{{{\left\| {x_k^i - x_{k - 1}^i} \right\|}_2^2}}}{{{L_f}{{\left\| {x_k^i - x_{k - 1}^i} \right\|}_2^2}}}.\\
 =&\frac{1}{{{mL_f}}}.
 \end{aligned}
\end{equation}
Next, we derive the bounds on BB step-size ${\alpha _k^i}$ following (\ref{E7}). According to Lipschitz continuity, it holds that
\begin{equation}
{L_f}{\left( {x_k^i - x_{k - 1}^i} \right)^{\top } }\left( {\nabla {f_i}\left( {x_k^i} \right) - \nabla {f_i}\left( {x_{k - 1}^i} \right)} \right) \ge {\left\| {\nabla {f_i}\left( {x_k^i} \right) - \nabla {f_i}\left( {x_{k - 1}^i} \right)} \right\|_2^2}.
\end{equation}
Then, step-size ${\alpha _k^i}$ is lower bounded by $1/({m L_f})$. Clearly, the upper bound on step-size ${\alpha _k^i}$ is no greater than $1/({m \mu})$. The above proofs show that the range of step-size ${\alpha _k^i}$ computed by (\ref{E6}) is longer than the range of step-size ${\alpha _k^i}$ computed by (\ref{E7}). Recalling the definitions of ${\alpha _{\max }}$ and ${{\bar \alpha }_{\max }}$, we can directly obtain that
\begin{equation}\label{Em}
\frac{1}{{m{L_f}}} \le {\alpha _{\max }},{{\bar \alpha }_{\max }} \le \frac{1}{{m\mu }}.
\end{equation}
The proof is completed.
\end{proof}
\end{lemma}
\begin{remark}\label{R6}
It is worth mentioning that owing to the availability of larger step-sizes, ADBB is competitive than the state-of-art algorithms \cite{Nedic2017b,Xi2018a,Xi2018d,Xin2018,Xin2019d,Xin2019,Li2019c,Li2018,Li2019b} that can be applied into unbalanced directed networks. Especially, the largest step-sizes employed by \cite{Xi2018a,Xi2017b,Xi2018d,Xin2019} do not exceed $1/(m{L_f})$ which is exactly the lower bound of the step-size employed in ADBB.
\end{remark}

\begin{lemma}\label{L2}
Supposing that Assumption \ref{A1} holds, recalling the definitions of augmented weight matrix ${\mathcal{A}} = A \otimes {I_n}$ and ${{y_\infty }}$, for $\forall x \in {\mathbb{R}^{mn}}$, there holds
\begin{equation}\label{E12}
\left\| {{\mathcal{A}^H}x - {y_\infty }x} \right\|_{\pi} \le {\sigma ^H}\left\| {x - {y_\infty }x} \right\|_{\pi},
\end{equation}
where $0 < \sigma  := \left\| {\mathcal{A} - {y_\infty }} \right\|_{\pi} < 1$ serves as the mixing rate of the directed network.
\begin{proof}
Since $A$ is primitive and row-stochastic, it holds that ${y_\infty } = \left( {{1_m}\pi^{\top }} \right) \otimes {I_n}$ from Perron-Frobenius theorem \cite{Xin2019}. Then, in terms of (\ref{E10.2}), we get that ${\mathcal{A}^H}{y_\infty } = {y_\infty }$ and ${y_\infty }{y_\infty } = {y_\infty }$, which yields that
\begin{equation}
{\mathcal{A}^H}x - {y_\infty }x = \left( {{\mathcal{A}^H} - {y_\infty }} \right)\left( {x - {y_\infty }x} \right).
\end{equation}
We know that $\rho \left( {\mathcal{A} - {y_\infty }} \right) < 1$ and
\begin{equation}
\begin{aligned}
{\mathcal{A}^H} - {y_\infty } =& \left( {{\mathcal{A}^{H - 1}} - {y_\infty }} \right)\left( {\mathcal{A} - {y_\infty }} \right)\\
=&{\left( {\mathcal{A} - {y_\infty }} \right)^H}.
\end{aligned}
\end{equation}
Then, it is easy to obtain that
$\left\| {{\mathcal{A}^H} - {y_\infty }} \right\|_{\pi} < 1$.
The proof follows via setting $\sigma  = \left\| {\mathcal{A} - {y_\infty }} \right\|_{\pi}$.
\end{proof}
\end{lemma}

\begin{lemma}\label{L3}
Suppose that Assumption \ref{A1} holds. Recalling the definitions of ${y_k}$ and ${{y_\infty }}$, it holds
\begin{equation}\label{E13}
\left\| {{y_k} - {y_\infty }} \right\|_2 \le {\sigma ^{kH}},\forall k \ge 0,
\end{equation}
where $0 < {\sigma ^H} < 1$ is a constant.\newline
\begin{proof}
In terms of (\ref{E10.2}), we obtain that ${y_k} = {A^{kH}} \otimes {I_n} = {\mathcal{A}^{kH}}$.
Similar with the proof in Lemma \ref{L2}, it holds that
\begin{equation}
\left\| {{y_k} - {y_\infty }} \right\|_2 = \left\| {{{\left( {{\mathcal{A}^H} - {y_\infty }} \right)}^k}} \right\|_2 \le {\sigma ^{kH}}, \forall k \ge 0.
\end{equation}
The proof is completed.
\end{proof}
\end{lemma}

\begin{lemma}\label{L4}
Suppose that Assumption \ref{A1} holds. For $\forall k \ge 0$, the following inequalities hold
\begin{equation}\label{E14}
\left\| {\hat y_k^{ - 1} - \hat y_\infty ^{ - 1}} \right\|_2 \le \sqrt m {\sigma ^{kH}}{Y^2},
\end{equation}
\begin{equation}\label{E16}
\left\| {\hat y_{k + 1}^{ - 1} - \hat y_k^{ - 1}} \right\|_2 \le 2\sqrt m {\sigma ^{kH}}{Y^2}.
\end{equation}
\begin{proof}
We first give the proof of (\ref{E14}) as follows:
\begin{equation}
\begin{aligned}
  \left\| {\hat y_k^{ - 1} - \hat y_\infty ^{ - 1}} \right\|_2 =& \left\| {\hat y_k^{ - 1}\left( {{{\hat y}_\infty } - {{\hat y}_k}} \right)\hat y_\infty ^{ - 1}} \right\|_2\\
  \le& \left\| {\hat y_k^{ - 1}} \right\|_2\left\| {{\rm{diag}}\left\{ {{y_k} - {y_\infty }} \right\}} \right\|_2\left\| {\hat y_\infty ^{ - 1}} \right\|_2\\
   \le& \sqrt m {\sigma ^{kH}}{Y^2}.
  \end{aligned}
\end{equation}
The proof of (\ref{E16}) directly follows from the proof of (\ref{E14}), i.e.,
\begin{equation}
\begin{aligned}
  \left\| {\hat y_{k + 1}^{ - 1} - \hat y_k^{ - 1}} \right\|_2 =& \left\| {\hat y_{k + 1}^{ - 1} - \hat y_\infty ^{ - 1} + \hat y_\infty ^{ - 1} - \hat y_k^{ - 1}} \right\|_2\\
  \le& \left\| {\hat y_{k + 1}^{ - 1} - \hat y_\infty ^{ - 1}} \right\|_2 + \left\| {\hat y_\infty ^{ - 1} - \hat y_k^{ - 1}} \right\|_2\\
  \le& {\text{2}}\sqrt m {\sigma ^{kH}}{Y^2},
  \end{aligned}
\end{equation}
where the first inequality utilizes the triangle inequality.
\end{proof}
\end{lemma}

\begin{lemma}\label{L5}
Suppose that Assumption \ref{A1} holds. For $\forall k \ge 0$, the following equation holds
\begin{equation}\label{E17}
{y_\infty }{z_k} = {y_\infty }\hat y_k^{ - 1}\nabla F\left( {{x_k}} \right).
\end{equation}
\begin{proof}
Recalling that ${y_\infty }{\mathcal{A}^H} = {y_\infty }$, from (\ref{E10.3}), we recursively update
\begin{equation}
\begin{aligned}
  {y_\infty }{z_k} =& {y_\infty }{z_{k - 1}} + {y_\infty }\hat y_k^{ - 1}\nabla F\left( {{x_k}} \right) - {y_\infty }\hat y_{k - 1}^{ - 1}\nabla F\left( {{x_{k - 1}}} \right)\\
  =&{y_\infty }{z_{k - 2}} + {y_\infty }\hat y_{k - 1}^{ - 1}\nabla F\left( {{x_{k - 1}}} \right) - {y_\infty }\hat y_{k - 2}^{ - 1}\nabla F\left( {{x_{k - 2}}} \right)\\
  &+ {y_\infty }\hat y_k^{ - 1}\nabla F\left( {{x_k}} \right) - {y_\infty }\hat y_{k - 1}^{ - 1}\nabla F\left( {{x_{k - 1}}} \right)\\
   =&{y_\infty }\hat y_k^{ - 1}\nabla F\left( {{x_k}} \right),
   \end{aligned}
\end{equation}
where the last equality follows from the initial conditions that ${\hat y_0} = {I_{mn}}$ and ${z_0} = \nabla F\left( {{x_0}} \right)$.
\end{proof}
\end{lemma}
\begin{lemma}\label{L6}(\cite{Bubeck2015})
Let the global objective function $f\left( x \right)$ be $\mu$-strongly convex and $L_f$-Lipschitz continuous, where ${0 < \mu  \le {L_f}}$. Then, for $\forall x \in {\mathbb{R}^{m}}$ and $0 < \theta  < \frac{2}{{{L_f}}}$, we have
\begin{equation}\label{E18}
\left\| {x - \theta  \nabla f\left( x \right) - {x^*}} \right\|_2 \le \eta \left\| {x - {x^*}} \right\|_2,
\end{equation}
where $\eta  = \max \left\{ {\left| {1 - \theta \mu } \right|,\left| {1 - \theta {L_f}} \right|} \right\}$.
\end{lemma}

\begin{lemma}\label{L10}
Suppose that Assumption \ref{A1} holds. For $\forall k \ge 0$, the following inequality holds
\begin{equation}\label{E22}
\begin{aligned}
\left\| {{z_k}} \right\|_2 \le& {{\bar \pi }^{0.5}}m{L_f}\left\| {{x_k} - {y_\infty }{x_k}} \right\|_{\pi} + m{L_f}\left\| {{y_\infty }{x_k} - {1_m} \otimes {{\tilde x}^*}} \right\|_2\\
&+ {{\underline{\pi } }^{ - 0.5}} \left\| {{z_k} - {y_\infty }{z_k}} \right\|_{\pi} + \sqrt m \hat Y{Y^2}{\sigma ^{kH}}\left\| {\nabla F\left( {{x_k}} \right)} \right\|_2.
\end{aligned}
\end{equation}
\begin{proof}
Recalling ${y_\infty }\hat y_\infty ^{ - 1} = \left( {{1_m}{{1}}_m^{\top }} \right) \otimes {I_n}$ and ${y_\infty }{z_k} = {y_\infty }\hat y_k^{ - 1}\nabla F\left( {{x_k}} \right)$ from Lemma \ref{L5}, we obtain that
\begin{equation*}
\begin{aligned}
\left\| {{z_k}} \right\|_2 \le& \left\| {{z_k} - {y_\infty }{z_k}} \right\|_2 + \left\| {{y_\infty }{z_k}} \right\|_2\\
\le& {{\underline{\pi } }^{ - 0.5}}\left\| {{z_k}  -  {y_\infty }{z_k}} \right\|_{\pi}  + \left\| {{y_\infty }\hat y_k^{  -  1}\nabla F\left( {{x_k}} \right)  -  {y_\infty }\hat y_\infty ^{ - 1}\nabla F\left( {{x_k}} \right)} \right\|_2\\
&+ \left\| {{y_\infty }\hat y_\infty ^{ - 1}\nabla F\left( {{x_k}} \right) - \left( {{1_m}{{1}}_m^{\top }} \right) \otimes {I_n}\nabla F\left( {{{ x}^*}} \right)} \right\|_2, \\
\le& {{\underline{\pi } }^{ - 0.5}}\left\| {{z_k}  -  {y_\infty }{z_k}} \right\|_{\pi}  + {\left\| {{y_\infty }\hat y_k^{ - 1} - {y_\infty }\hat y_\infty ^{ - 1}} \right\|_2}{\left\| {\nabla F\left( {{x_k}} \right)} \right\|_2}\\
&+ m{\left\| {\nabla F\left( {{x_k}} \right) - \nabla F\left( {{x^*}} \right)} \right\|_2}, \\
\end{aligned}
\end{equation*}
where the second inequality utilizes the fact that $\left( {1_m^{\top } \otimes {I_n}} \right)\nabla F\left( {{x^*}} \right) = 0$. Applying Lipschitz continuity and (\ref{E14}) completes the proof.
\end{proof}
\end{lemma}

To proceed, the contraction relationship of ADBB is established by deriving upper bounds on the sequel quantities: i) the consensus error: $\left\| {{x_{k + 1}} - {y_\infty }{x_{k + 1}}} \right\|_{ \pi}$; ii) the optimality gap: $\left\| {{y_\infty }{x_{k + 1}} - {1_m} \otimes {{\tilde x}^*}} \right\|_2$; iii)  the gradient tracking error: $\left\| {{z_{k + 1}} - {y_\infty }{z_{k + 1}}} \right\|_{ \pi}$.

\subsection{Contraction Relationship}
\begin{lemma}\label{L7}
Suppose that Assumption \ref{A1} holds. For $\forall k \ge 0$, the following inequality holds
\begin{equation}\label{E19}
\begin{aligned}
&\left\| {{x_{k + 1}} - {y_\infty }{x_{k + 1}}} \right\|_{ \pi}\\
\le & \left( {1 + {\alpha _{\max }}\bar \pi m{L_f}} \right){\sigma ^H}{\left\| {{x_k} - {y_\infty }{x_k}} \right\|_\pi } + {\alpha _{\max }}{{\bar \pi }^{0.5}}m{L_f}{\sigma ^H}{\left\| {{y_\infty }{x_k} - {1_m} \otimes {{\tilde x}^*}} \right\|_2}\\
&+ {\alpha _{\max }}{\vartheta ^{0.5}}{\sigma ^H} {\left\| {{z_k} - {y_\infty }{z_k}} \right\|_\pi } + {\alpha _{\max }}\sqrt {m\bar \pi } \hat Y{Y^2}{\sigma ^{\left( {k + 1} \right)H}}{\left\| {\nabla F\left( {{x_k}} \right)} \right\|_2}.
\end{aligned}
\end{equation}
\begin{proof}
According to (\ref{E10.1}) and ${y_\infty }{\mathcal{A}^H} = {y_\infty }$, it holds that
\begin{equation}\label{E19+}
\begin{aligned}
&\left\| {{x_{k + 1}} - {y_\infty }{x_{k + 1}}} \right\|_{ \pi} \\
=& \left\| {{\mathcal{A}^H}\left( {{x_k} - D_k^\alpha {z_k}} \right) - {y_\infty }\left( {{x_k} - D_k^\alpha {z_k}} \right)} \right\|_{ \pi}\\
=&\left\| {\left( {{\mathcal{A}^H} - {y_\infty }} \right)\left( {{x_k} - {y_\infty }{x_k}} \right) + \left( {{\mathcal{A}^H} - {y_\infty }} \right)D_k^\alpha {z_k}} \right\|_{ \pi}\\
\le& {\sigma ^H}\left\| {{x_k} - {y_\infty }{x_k}} \right\|_{ \pi} + {\alpha _{\max }}{{\bar \pi }^{0.5}}{\sigma ^H}\left\| {{z_k}} \right\|_2,
\end{aligned}
\end{equation}
where the second equality follows the relationship $\left( {{\mathcal{A}^H} - {y_\infty }} \right){y_\infty } = 0$ and the inequality applies Lemma \ref{L2}. The proof is finished through substituting (\ref{E22}) into (\ref{E19+}).
\end{proof}
\end{lemma}

\begin{lemma}\label{L8}
Suppose that Assumptions \ref{A1}-\ref{A3} hold. If $\pi^{\top }{\alpha_k}  < 2/\left( {m{L_f}} \right)$, for $\forall k \ge 0$, then the following inequality holds
\begin{equation}\label{E20}
\begin{aligned}
&\left\| {{y_\infty }{x_{k + 1}} - {1_m} \otimes {{\tilde x}^*}} \right\|_2\\
\le& {\alpha _{\max }}{{\underline{\pi } }^{ - 0.5}}{mL_f}\left\| {{x_k} - {y_\infty }{x_k}} \right\|_{ \pi} + \lambda \left\| {{y_\infty }{x_k} - {1_m} \otimes {{\tilde x}^*}} \right\|_2\\
& + {\alpha _{\max }}{{\underline{\pi } }^{ - 0.5}}\hat Y\left\| {{z_k} - {y_\infty }{z_k}} \right\|_{ \pi} + {\alpha _{\max }}\sqrt m \hat Y{Y^2}{\sigma ^{kH}}\left\| {\nabla F\left( {{x_k}} \right)} \right\|_2,
\end{aligned}
\end{equation}
where $\lambda  = \max \left\{ {\left| {1 - m\pi^{\top }{\alpha _k}\mu } \right|,\left| {1 - m\pi^{\top }{\alpha _k}{L_f}} \right|} \right\}$. \newline
\begin{proof}
Recalling ${y_\infty }\mathcal{A}{^H} = {y_\infty }$, we obtain that
\begin{equation}\label{E20-1}
\begin{aligned}
&\left\| {{y_\infty }{x_{k + 1}} - {1_m} \otimes {{\tilde x}^*}} \right\|_2 \\
 =&  \left\| {{y_\infty }\left( {{\mathcal{A}^H}\left( {{x_k} - D_k^\alpha {z_k}} \right) + \left( {D_k^\alpha  - D_k^\alpha } \right){y_\infty }{z_k}} \right) - {1_m} \otimes {{\tilde x}^*}} \right\|_2\\
\le& \left\| {{y_\infty }{x_k}  -  {y_\infty }D_k^\alpha {y_\infty }{z_k}  -  {1_m} \otimes {{\tilde x}^*}} \right\|_2 +  {\alpha _{\max }}\hat Y\left\| {{z_k} - {y_\infty }{z_k}} \right\|_2\\
\le& \left\| {{y_\infty }{x_k}  -  {y_\infty }D_k^\alpha {y_\infty }{z_k}  -  {1_m} \otimes {{\tilde x}^*}} \right\|_2 +  {\alpha _{\max }}\hat Y{{\underline{\pi } }^{ - 0.5}}\left\| {{z_k} - {y_\infty }{z_k}} \right\|_{\pi}.
\end{aligned}
\end{equation}
We know that
\begin{equation}
\begin{aligned}
{y_\infty }D_k^\alpha {y_\infty } =& \left( {\left( {{1_m}\pi^{\top }} \right) \otimes  {I_n}} \right)\left( {{\rm{diag}}\left\{ {{\alpha _k}} \right\}  \otimes  {I_n}} \right)\left( {\left( {{1_m}\pi^{\top }} \right)  \otimes  {I_n}} \right)\\
 =& \left( {\pi^{\top }{\alpha _k}} \right){y_\infty }.
\end{aligned}
\end{equation}
We continue to bound $\left\| {{y_\infty }{x_k} - {y_\infty }D_k^\alpha {y_\infty }{z_k} - {1_m} \otimes {{\tilde x}^*}} \right\|_2$ as follows:
\begin{equation}\label{E20-2}
\begin{aligned}
&\left\| {{y_\infty }{x_k} - {y_\infty }D_k^\alpha {y_\infty }{z_k} - {1_m} \otimes {{\tilde x}^*}} \right\|_2\\
 \le& \sqrt m \left\| {\left( {\pi^\top   \otimes {I_n}} \right){x_k}  -  {{\tilde x}^*}  -  m\left( {\pi^\top  {\alpha _k}} \right)\nabla f\left( {\left( {\pi^\top   \otimes {I_n}} \right){x_k}} \right)} \right\|_2\\
 &+\left\| {m\left( {\pi^\top  {\alpha _k}} \right)\left( {{1_m} \otimes {I_n}} \right)\nabla f\left( {\left( {\pi^\top   \otimes {I_n}} \right){x_k}} \right)} { - \left( {\pi^\top  {\alpha _k}} \right){y_\infty }{z_k}} \right\|_2\\
  =&{\Delta_1} + {\Delta_2}.
\end{aligned}
\end{equation}
If ${\pi ^ \top }{\alpha _k} < 2/\left( {m{L_f}} \right)$, then using Lemma \ref{L6}, it holds that
\begin{equation}\label{E20-3}
\begin{aligned}
{\Delta_1} \le \lambda \left\| {{y_\infty }{x_k} - {1_m} \otimes {{\tilde x}^*}} \right\|_2,
\end{aligned}
\end{equation}
where $\lambda  = \max \left\{ {\left| {1 - m\pi^{\top }{\alpha _k}\mu } \right|,\left| {1 - m\pi^{\top }{\alpha _k}{L_f}} \right|} \right\}$. Next, we aim to bound ${\Delta_2}$.
\begin{equation}\label{E20-4}
\begin{aligned}
{\Delta_2} =& \left( {\pi^{\top }{\alpha _k}} \right)\| m{\left( {{1_m} \otimes {I_n}} \right)\nabla f\left( {\left( {\pi^{\top } \otimes {I_n}} \right){x_k}} \right) - {y_\infty }{z_k}} \|_2\\
\le& {\alpha _{\max }}\| m{\left( {{1_m} \otimes {I_n}} \right)\nabla f\left( {\left( {\pi^\top   \otimes {I_n}} \right){x_k}} \right)}{ - \left( {{1_m} \otimes {I_n}} \right)\left( {{1}_m^\top   \otimes {I_n}} \right)\nabla F\left( {{x_k}} \right)} \|_2\\
&+ {\alpha _{\max }}\| {\left( {{1_m} \otimes {I_n}} \right)\left( {{{1}}_m^{\top } \otimes {I_n}} \right)\nabla F\left( {{x_k}} \right) - {y_\infty }\hat y_k^{ - 1}\nabla F\left( {{x_k}} \right)} \|_2\\
\le& {\alpha _{\max }}{mL_f}\left\| {{x_k} - {y_\infty }{x_k}} \right\|_2 + {\alpha _{\max }} \sqrt m \hat Y{Y^2}{\sigma ^{kH}}\left\| {\nabla F\left( {{x_k}} \right)} \right\|_2\\
\le& {\alpha _{\max }}{{\underline{\pi } }^{ - 0.5}}{mL_f}\left\| {{x_k} - {y_\infty }{x_k}} \right\|_{\pi} + {\alpha _{\max }} \sqrt m \hat Y{Y^2}{\sigma ^{kH}}\left\| {\nabla F\left( {{x_k}} \right)} \right\|_2,
\end{aligned}
\end{equation}
where the first inequality uses Lemma \ref{L5} and the last inequality applies Lipschitz continuity and (\ref{E14}). Plugging (\ref{E20-3}) and (\ref{E20-4}) into (\ref{E20-2}) yields
\begin{equation}\label{E20-5}
\begin{aligned}
&\left\| {{y_\infty }{x_k} - {y_\infty }D_k^\alpha {y_\infty }{z_k} - {1_m} \otimes {{\tilde x}^*}} \right\|_2\\
\le & {\alpha _{\max }}{{\underline{\pi } }^{ - 0.5}}{mL_f}\left\| {{x_k} - {y_\infty }{x_k}} \right\|_{\pi} + \lambda \left\| {{y_\infty }{x_k} - {1_m} \otimes {{\tilde x}^*}} \right\|_2\\
& +  {\alpha _{\max }} \sqrt m \hat Y{Y^2}{\sigma ^{kH}}\left\| {\nabla F\left( {{x_k}} \right)} \right\|_2.\\
\end{aligned}
\end{equation}
The proof is completed by substituting (\ref{E20-5}) into (\ref{E20-1}).
\end{proof}
\end{lemma}

\begin{lemma}\label{L9}
Suppose that Assumptions \ref{A1}-\ref{A2} hold. For $\forall k \ge 0$, the following inequality holds
\begin{equation}\label{E21}
\begin{aligned}
\left\| {{z_{k + 1}} - {y_\infty }{z_{k + 1}}} \right\|_{\pi}\le& \left( {{\vartheta ^{0.5}}Y{L_f}{p_1} + {\alpha _{\max }}{{\bar \pi }^{0.5}}{mY\hat YL_f^2}} \right){\sigma ^H}{\left\| {{x_k} - {y_\infty }{x_k}} \right\|_\pi }\\
&+ {\alpha _{\max }}mY\hat YL_f^2{\sigma ^H}{\left\| {{y_\infty }{x_k} - {1_m} \otimes {{\tilde x}^*}} \right\|_2} \\
&+ \left( {1 + {\alpha _{\max }}{{\underline{\pi } }^{ - 0.5}}Y\hat Y{L_f}} \right){\sigma ^H}{\left\| {{z_k} - {y_\infty }{z_k}} \right\|_\pi }\\
&+ \left( {2\sqrt {m{{\underline{\pi } }^{ - 1}}}   + {\alpha _{\max }}\sqrt m {L_f}Y{{\hat Y}^2}} \right){Y^2}{\sigma ^{\left( {k + 1} \right)H}}{\left\| {\nabla F\left( {{x_k}} \right)} \right\|_2}
\end{aligned}
\end{equation}
\begin{proof}
According to Lemma \ref{L2} and (\ref{E10.3}), we obtain that
\begin{equation}\label{E21-1}
\begin{aligned}
&\left\| {{z_{k + 1}} - {y_\infty }{z_{k + 1}}} \right\|_{\pi}\\
=& \left\| {{\mathcal{A}^H}\left( {{z_k} + \hat y_{k + 1}^{ - 1}\nabla F\left( {{x_{k + 1}}} \right) - \hat y_k^{ - 1}\nabla F\left( {{x_k}} \right)} \right) - {y_\infty }{z_{k + 1}}} \right\|_{\pi}\\
 \le& \left\| {{\mathcal{A}^H}\left( {y_{k  +  1}^{  -  1}\nabla F\left( {{x_{k  +  1}}} \right)  -  \hat y_k^{  -  1}\nabla F\left( {{x_k}} \right)} \right)  -  \left( {{y_\infty }{z_{k  +  1}}  -  {y_\infty }{z_k}} \right)} \right\|_{\pi}\\
 &+{\sigma ^H}\left\| {{z_k} - {y_\infty }{z_k}} \right\|_{\pi}.
 \end{aligned}
\end{equation}
Recall that ${y_\infty }{z_k} = {y_\infty }\hat y_k^{ - 1}\nabla F\left( {{x_k}} \right)$ from Lemma \ref{L5}. Therefore
\begin{equation}\label{E21-2}
\begin{aligned}
&\left\| {{\mathcal{A}^H}\left( {y_{k + 1}^{ - 1}\nabla F\left( {{x_{k +  1}}} \right) - \hat y_k^{  -  1}\nabla F\left( {{x_k}} \right)} \right)  -  \left( {{y_\infty }{z_{k  +  1}} - {y_\infty }{z_k}} \right)} \right\|_{\pi}\\
=& \left\| {\left( {{\mathcal{A}^H} - {y_\infty }} \right)\left( {y_{k + 1}^{ - 1}\nabla F\left( {{x_{k + 1}}} \right) - \hat y_k^{ - 1}\nabla F\left( {{x_k}} \right)} \right)} \right\|_{\pi}\\
\le& {\sigma ^H}\left\| {y_{k + 1}^{ - 1}\nabla F\left( {{x_{k + 1}}} \right) - y_{k + 1}^{ - 1}\nabla F\left( {{x_k}} \right)} \right\|_{\pi} + {\sigma ^H}\left\| {y_{k + 1}^{ - 1}\nabla F\left( {{x_k}} \right) - \hat y_k^{ - 1}\nabla F\left( {{x_k}} \right)} \right\|_{\pi}\\
\le& {{\underline{\pi } }^{ - 0.5}}Y{L_f}{\sigma ^H}\left\| {{x_{k + 1}} - {x_k}} \right\|_2 + 2{{\underline{\pi } }^{ - 0.5}}\sqrt m{Y^2}{\sigma ^{(k + 1)H}}\left\| {\nabla F\left( {{x_k}} \right)} \right\|_2.
\end{aligned}
\end{equation}
We next bound $\left\| {{x_{k + 1}} - {x_k}} \right\|$,
\begin{equation}\label{E21-3}
\begin{aligned}
&\left\| {{x_{k + 1}} - {x_k}} \right\|_2 \\
=& \left\| {{\mathcal{A}^H}\left( {{x_k} - D_k^\alpha {z_k}} \right) - {x_k}} \right\|_2\\
\le& \left\| {\left( {{I_{mn}} - {\mathcal{A}^H}} \right)\left( {{x_k} - {y_\infty }{x_k}} \right)} \right\|_2 + {\alpha _{\max }}\left\| {{\mathcal{A}^H}} \right\|_2\left\| {{z_k}} \right\|_2\\
\le& {{\bar \pi }^{0.5}}{p_1}\left\| {{x_k} - {y_\infty }{x_k}} \right\|_{\pi} + {\alpha _{\max }}\hat Y\left\| {{z_k}} \right\|_2,
\end{aligned}
\end{equation}
where the first inequality utilizes the fact that $\left( {{I_{mn}} - {\mathcal{A}^H}} \right){y_\infty } = 0$. Combining (\ref{E21-2}) and (\ref{E21-3}) gives
\begin{equation}\label{E21-4}
\begin{aligned}
&\left\| {{\mathcal{A}^H}\left( {y_{k + 1}^{ - 1}\nabla F\left( {{x_{k +  1}}} \right) - \hat y_k^{  -  1}\nabla F\left( {{x_k}} \right)} \right)  -  \left( {{y_\infty }{z_{k  +  1}} - {y_\infty }{z_k}} \right)} \right\|_{\pi}\\
\le& {\vartheta ^{0.5}}{p_1}Y{L_f}{\sigma ^H}{\left\| {{x_k} - {y_\infty }{x_k}} \right\|_\pi } \!+\! 2{\sqrt {m\underline{\pi }^{ - 1} }} {Y^2}{\sigma ^{(k + 1)H}}{\left\| {\nabla F\left( {{x_k}} \right)} \right\|_2} \!+\! {\alpha _{\max }}\hat Y{\left\| {{z_k}} \right\|_2}.
\end{aligned}
\end{equation}
Substituting (\ref{E22}) into (\ref{E21-1}) gives
\begin{equation}\label{E21-5}
\begin{aligned}
&\left\| {{z_{k + 1}} - {y_\infty }{z_{k + 1}}} \right\|_{\pi}\\
\le& {\vartheta ^{0.5}}{p_1}Y{L_f}{\sigma ^H}{\left\| {{x_k} - {y_\infty }{x_k}} \right\|_\pi } +{\sigma ^H}\left\| {{z_k} - {y_\infty }{z_k}} \right\|_{\pi} \\
&+ 2{\sqrt {m\underline{\pi }^{ - 1} }} {Y^2}{\sigma ^{(k + 1)H}}{\left\| {\nabla F\left( {{x_k}} \right)} \right\|_2}+ {\alpha _{\max }}\hat Y{\left\| {{z_k}} \right\|_2}.
\end{aligned}
\end{equation}
The proof is ended via substituting (\ref{E22}) into (\ref{E21-5}).
\end{proof}
\end{lemma}

\begin{lemma}\label{L11}(\cite[Lemma 4]{Xin2019d})
Let $X \in {\mathbb{R}^{n \times n}}$ be a non-negative matrix and $x \in {\mathbb{R}^n}$ be a positive vector. If $Xx < \omega x$, then $\rho \left( X \right) < \omega$.
\end{lemma}

Before presenting the main results, we make some definitions for subsequent proofs as follows:
\[\begin{array}{*{20}{c}}
  {\begin{array}{*{20}{c}}
  {\begin{array}{*{20}{c}}
  {{t_k} = \left[ {\begin{array}{*{20}{c}}
  {\left\| {{x_k} - {y_\infty }{x_k}} \right\|_{\pi}} \\
  {\left\| {{y_\infty }{x_k} - {1_m} \otimes {{\tilde x}^*}} \right\|_2} \\
  {\left\| {{z_k} - {y_\infty }{z_k}} \right\|_{\pi}}
\end{array}} \right],\;\;\;\,\,{g_k} = \left[ {\begin{array}{*{20}{c}}
  {\left\| {\nabla F\left( {{x_k}} \right)} \right\|_2} \\
  0 \\
  0
\end{array}} \right],}
\end{array}}
\end{array}}
\end{array}\]
\[{{M_k} = \left[ {\begin{array}{*{20}{c}}
  {{\sigma ^H} + {\alpha _{\max }}{\sigma ^H}{w_1}}&{{\alpha _{\max }}{\sigma ^H}{w_1}}&{{\alpha _{\max }}{\sigma ^H}} \\
  {{\alpha _{\max }}{w_2}}&\lambda &{{\alpha _{\max }}{w_3}} \\
  {{\sigma ^H}{w_4} + {\alpha _{\max }}{\sigma ^H}{w_5}}&{{\alpha _{\max }}{\sigma ^H}{w_5}}&{{\sigma ^H} + {\alpha _{\max }}{\sigma ^H}{w_6}}
\end{array}} \right]},\]
\[{G_k} = \left[ {\begin{array}{*{20}{c}}
  {{\alpha _{\max }}\sqrt {m\bar \pi } \hat Y{Y^2}{\sigma ^H}}&0&0 \\
  {{\alpha _{\max }}\sqrt m {Y^2}\hat Y}&0&0 \\
  {2\sqrt {m{{\underline{\pi } }^{ - 1}}}  {Y^2}{\sigma ^H} + {\alpha _{\max }}\sqrt m {{\hat Y}^2}{Y^3}{L_f}{\sigma ^H}}&0&0
\end{array}} \right]{\sigma ^{kH}},\]
where ${w_1} =  {{\bar \pi }^{0.5}} m{L_f}$, ${w_2} = {{\underline{\pi } }^{ - 0.5}}{mL_f}$, ${w_3} = {{\underline{\pi } }^{ - 0.5}}\hat Y$, ${w_4} = {{\vartheta ^{0.5}}}Y{L_f}{p_1}$, ${w_5} = mY\hat YL_f^2$, ${w_6} = {{{\underline{\pi } }^{ - 0.5}}}Y\hat Y{L_f}$. Note that ${G_k}$ decays linearly over $k$, since $0 < \sigma  < 1$.
\subsection{Main Results}
\begin{theorem}\label{T1}
Suppose that Assumptions \ref{A1}-\ref{A3} hold. If $\pi^{\top }{\alpha _k} < 2/\left( {m{L_f}} \right)$ for $\forall k \ge 0$, one can establish a linear time-variant inequality as follows:
\begin{equation}\label{E23}
\begin{aligned}
{t_{k + 1}} \le M_k {t_k} + {G_k}{g_k},\forall k \ge 0,
\end{aligned}
\end{equation}
where ${t_k},{g_k} \in {\mathbb{R}^3}$, and $M_k ,{G_k} \in {\mathbb{R}^{3 \times 3}}$ are defined before Theorem \ref{T1}. Furthermore, when the largest step-size ${\alpha _{\max }}$ satisfies (\ref{Em}) with proper setting of the multi-consensus inner loop number $H$, there holds that $\rho \left( {M_k} \right) < 1$. \newline
\begin{proof}
Recall that $\lambda  = \max \left\{ {\left| {1 - m\pi^{\top }{\alpha _k}\mu } \right|,\left| {1 - m\pi^{\top }{\alpha _k}{L_f}} \right|} \right\}$. Clearly, we can obtain that $1/{mL_f} \le {\alpha _{\min }} \le \pi^{\top }{\alpha _k}$ from Lemma \ref{L1}. If $\pi^{\top }{\alpha _k} \le \left( {2/{mL_f}} \right) - \left( {\mu /{m}L_f^2} \right)$, it holds that $\lambda  \le 1 - \left( {\mu /{mL_f}} \right)$.
Hence, we get that $M_k  \le {M }$, with
\[M = \left[ {\begin{array}{*{20}{c}}
  {{\sigma ^H} + {\alpha _{\max }}{\sigma ^H}{w_1}}&{{\alpha _{\max }}{\sigma ^H}{w_1}}&{{\alpha _{\max }}{\sigma ^H}} \\
  {{\alpha _{\max }}{w_2}}&{1 - {w_7}}&{{\alpha _{\max }}{w_3}} \\
  {{\sigma ^H}{w_4} + {\alpha _{\max }}{\sigma ^H}{w_5}}&{{\alpha _{\max }}{\sigma ^H}{w_5}}&{{\sigma ^H} + {\alpha _{\max }}{\sigma ^H}{w_6}}
\end{array}} \right],\]
where ${w_7} = \mu /{L_f}$. Thus, $\rho \left( {M_k } \right) \le \rho \left( {{M }} \right)$. To proceed, we derive the lower bound on multi-consensus inner loop number $H$ and solve for a positive vector $c = {\left[ {{c_1},{c_2},{c_3}} \right]^{\top  }}$ from
\begin{equation}\label{E24}
{M }\left[ {\begin{array}{*{20}{c}}
  {{c_1}} \\
  {{c_2}} \\
  {{c_3}}
\end{array}} \right] < \left[ {\begin{array}{*{20}{c}}
  {{c_1}} \\
  {{c_2}} \\
  {{c_3}}
\end{array}} \right],
\end{equation}
which is equivalent to the following inequalities
\begin{equation}\label{E25}
\left\{ \begin{gathered}
  0 < {\alpha _{\max }}\left( {{\sigma ^H}{w_1}{c_1} + {\sigma ^H}{w_1}{c_2} + {\sigma ^H}{c_3}} \right) < \left( {1 - {\sigma ^H}} \right){c_1}, \hfill \\
  0 < {\alpha _{\max }}\left( {{w_2}{c_1} + {w_3}{c_3}} \right) < {w_7}{c_2}, \hfill \\
  0 < {\alpha _{\max }}\left( {{\sigma ^H}{w_5}{c_1} + {\sigma ^H}{w_5}{c_2} + {\sigma ^H}{w_6}{c_3}} \right) < \left( {1 - {\sigma ^H}} \right){c_3}- {\sigma ^H}{w_4}{c_1}. \hfill \\
\end{gathered}  \right.
\end{equation}
Since the right hand side of the last inequality in (\ref{E25}) has to be positive, it should hold that
\begin{equation}\label{E26}
0 < {c_1} < \frac{{\left( {1 - {\sigma ^H}} \right){c_3}}}{{{\sigma ^H}{w_4}}}.
\end{equation}
According to (\ref{E25}), we can obtain the upper bound on the largest step-size as follows:
\begin{equation}\label{E27}
0 \!<\! {\alpha _{\max }} \!<\! \min {\text{ }}\left\{ {\frac{{\left( {1 \!-\! {\sigma ^H}} \right){c_1}}}{{{\sigma ^H}\left( {{w_1}{c_1} \!+\! {w_1}{c_2} \!+\! {c_3}} \right)}},\frac{{{w_7}{c_2}}}{{{w_2}{c_1} \!+\! {w_3}{c_3}}},\frac{{\left( {1 \!-\! {\sigma ^H}} \right){c_3} \!-\! {\sigma ^H}{w_4}{c_1}}}{{{\sigma ^H}\left( {{w_5}{c_1} \!+\! {w_5}{c_2} \!+\! {w_6}{c_3}} \right)}}} \right\},
\end{equation}
where ${c_2},{c_3}$ are arbitrary positive constants and $c_1$ is chosen from (\ref{E26}). According to (\ref{Em}), we need
\begin{equation}\label{E28}
\frac{1}{{{mL_f}}} \le {\alpha _{\max }} \le \frac{1}{m\mu }.
\end{equation}
So as to guarantee that the range of largest step-size in (\ref{E27}) contains the range of the largest step-size in (\ref{E28}), the following inequality should hold
\begin{equation}\label{E29}
0 \!<\! \frac{1}{m\mu } \!<\! \min \left\{ {\frac{{\left( {1 \!-\! {\sigma ^H}} \right){c_1}}}{{{\sigma ^H}\left( {{w_1}{c_1} \!+\! {w_1}{c_2} \!+\! {c_3}} \right)}},\frac{{{w_7}{c_2}}}{{{w_2}{c_1} \!+\! {w_3}{c_3}}},\frac{{\left( {1 \!-\! {\sigma ^H}} \right){c_3} \!-\! {\sigma ^H}{w_4}{c_1}}}{{{\sigma ^H}\left( {{w_5}{c_1} \!+\! {w_5}{c_2} \!+\! {w_6}{c_3}} \right)}}} \right\},
\end{equation}
which yields
\begin{equation}\label{E30}
0 \!<\! {\sigma ^H} \!<\! \min \left\{ {\frac{{m\mu {c_1}}}{{\left( {{w_1} \!+\! m\mu } \right){c_1} \!+\! {w_1}{c_2} \!+\! {c_3}}},\frac{{m\mu {c_3}}}{{\left( {{w_5}{\!+\!}m\mu {w_4}} \right){c_1} \!+\! {w_5}{c_2} \!+\! \left( {{w_6} \!+\! m\mu } \right){c_3}}}} \right\},
\end{equation}
We define $\varpi  = \min \left\{ {\frac{{m\mu {c_1}}}{{\left( {{w_1} \!+\! m\mu } \right){c_1} \!+\! {w_1}{c_2} \!+\! {c_3}}},\frac{{m\mu {c_3}}}{{\left( {{w_5}{\!+\!}m\mu {w_4}} \right){c_1} \!+\! {w_5}{c_2} \!+\! \left( {{w_6} \!+\! m\mu } \right){c_3}}}} \right\}$. Then, it follows from (\ref{E30}) that
\begin{equation}\label{E31}
H \ge \left\lceil {\frac{{\ln \varpi }}{{\ln \sigma }}} \right\rceil  + \varphi \left( {\frac{{\ln \varpi }}{{\ln \sigma }}} \right),
\end{equation}
where $\varphi \left( x \right) = \left\{ \begin{gathered}
  1,{\text{for }}x \in {\mathbb{N}_ + } \hfill \\
  0,{\text{for }}x \notin {\mathbb{N}_ + } \hfill \\
\end{gathered}  \right.$ and ${\mathbb{N}_ + }$ denotes the set of positive integers. Thus, it holds that $\rho \left( {{M }} \right) < 1$. Finally, we get $\rho \left( {{M_k}} \right) < 1$ due to $\rho \left( {M_k } \right) \le \rho \left( {{M }} \right)$. The proof is completed.
\end{proof}
\end{theorem}
\begin{remark}\label{R7}
We emphasize that the BB step-size (see, (\ref{E6}) and (\ref{E7})) is automatically computed by each agent in the system and only depends on the local information, which does not rely on the network parameter, $\sigma$. We also acknowledge that multi-consensus inner loop number $H$ is lower bounded by the network parameter, $\sigma$, which cannot be computed locally. Therefore, we pick a sufficiently large multi-consensus number, $H$, to simultaneously guarantee the convergence of ADBB and ensure that ADBB can run in a fully distributed manner. This notion of sufficiently large only serves for the conservative theoretical analysis while  in practice we need to manually optimize the multi-consensus number, $H$, which is shown in Section \ref{section five}.
\end{remark}
\begin{lemma}\label{L12}(\cite[Theorem 2]{Xi2018a,Xin2019})
Suppose that Assumptions 1-3 hold.
Recall that $\rho \left( {M_k } \right) < 1$ according to Theorem \ref{T1} and ${G_k}$ decays linearly over $k$ since $0 < \sigma  < 1$. Then, the sequence ${\left\{ {{x_k}} \right\}_{k \ge 0}}$ generated by ADBB converges linearly to ${x^*} = {1_m} \otimes {{\tilde x}^*}$, i.e., for some positive constant $\omega  > 0$, it holds that
\begin{equation}\label{E32}
\left\| {{x_k} - {1_m} \otimes {{\tilde x}^*}} \right\| \le \omega {\left( {\max \left\{ {\rho \left( {{M_k}} \right),{\sigma ^H}} \right\} + \xi } \right)^k}, \forall k \ge 0,
\end{equation}
where $\xi $ is an arbitrary small constant such that $0 < \xi  + \max \left\{ {\rho \left( M_k \right),{\sigma ^H}} \right\} < 1$.
\end{lemma}
\section{Numerical Experiments}\label{section five}
In order to confirm correctness of the theoretical analysis and show the performance of ADBB. We leverage regularized logistic regression to compare ADBB with some other well-known distributed optimization algorithms through solving a binary classification problem. Especially, a network of $m$ agents cooperatively resolve a distributed logistic regression problem as follows:
\begin{equation}\label{E33}
\mathop {\min }\limits_{\tilde w \in {\mathbb{R}^n}} f\left( {\tilde w} \right) = \frac{1}{m}\sum\limits_{i = 1}^m {{f_i}\left( {\tilde w} \right)},
\end{equation}
where $\tilde w \in {\mathbb{R}^n}$ is the optimization variable to learn the separating hyperplane and each local objective function ${{f_i}\left( {\tilde w} \right)}$ is expressed as
\begin{equation}\label{E34}
 {f_i}\left( {\tilde w} \right) =  \frac{1}{{{q_i}}}\sum\limits_{j = 1}^{{q_i}} {\log \left( {1 + \exp \left( { - {b_{ij}}c_{ij}^{\top }\tilde w} \right)} \right)} + \frac{\beta }{{2m}}{\left\| {\tilde w} \right\|_2^2},
\end{equation}
where $q_i$ is the number of local samples distributed to agent $i$; ${c_{ij}} \in {\mathbb{R}^n}$ is the $j$-th training sample and ${b_{ij}} \in \left\{ { + 1, - 1} \right\}$ is the corresponding label, both of which are only accessed by agent $i$; $\beta$ is the regularized constant.
\begin{figure}[h!]
\begin{minipage}[h!]{0.5\linewidth}
{\includegraphics[height=5.5cm,width=6.5cm]{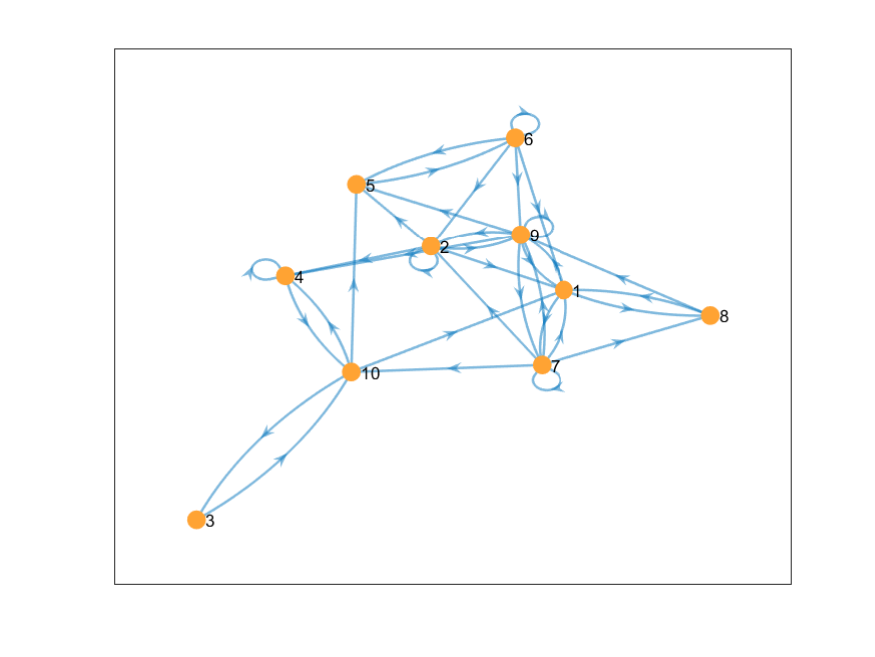}}
\caption{An unbalanced directed network with 10 agents.}
\label{Fig.1}
\end{minipage}
\begin{minipage}[h!]{0.5\linewidth}
{\includegraphics[height=5.5cm,width=6.5cm]{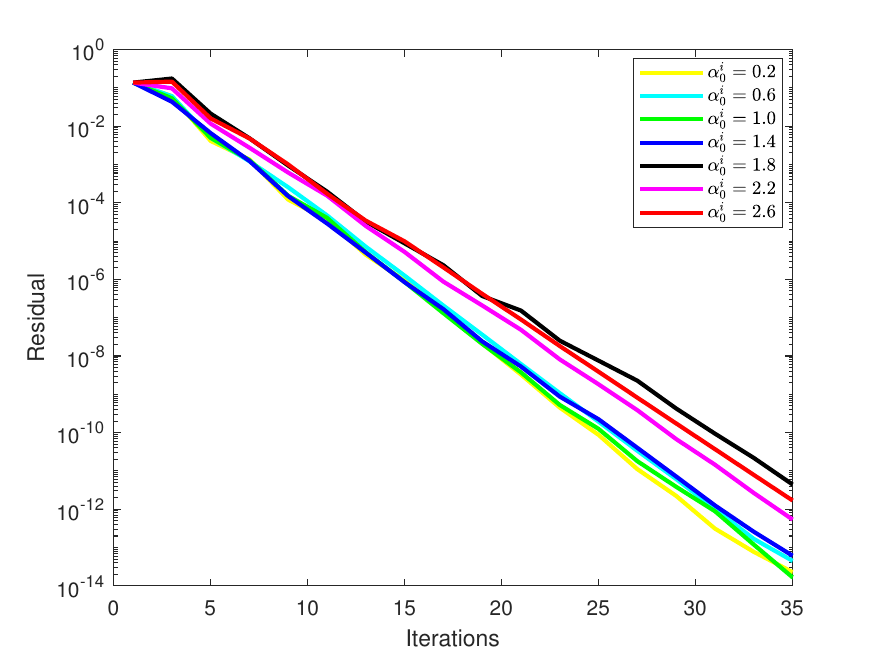}}
\caption{Performance of ADBB under different initial step-sizes $\alpha _0^i$.}
\label{Fig.2}
\end{minipage}
\end{figure}
Therefore, the optimal solution to (\ref{E33}) is represented by
\begin{equation}\label{E34+}
{{\tilde w}^*} = \mathop {\arg \min }\limits_{\tilde w \in {\mathbb{R}^n}} \left( {\frac{\beta }{2}{{\left\| {\tilde w} \right\|}_2^2} + \frac{1}{{{m}}}\sum\limits_{i = 1}^m \frac{1}{{{q_i}}}{\sum\limits_{j = 1}^{{q_i}} {\log \left( {1 + \exp \left( { - {b_{ij}}c_{ij}^{\top } \tilde w} \right)} \right)} } } \right).
\end{equation}
In the following experiments, we assume that the samples are distributed equally among the agents, i.e., ${q_i} = N/m$, $i \in \mathcal{V}$, where $N$ is the total number of data and $m$ represents the number of agents in the network. Then, we provide two case studies, in which the residual is defined as $\left( {1/m} \right)\sum\nolimits_{i = 1}^m {\left\| {w_k^i - {{\tilde w}^*}} \right\|_2}$. All simulations are carried out in MATLAB on a Lenovo laptop with 4.10 GHz, 4 Cores 8 Threads Intel i5 - 9300HF processor and 16GB memory.
\subsection*{Case Study 1}
In the first case, the effect of the multi-consensus inner loop number $H$ and initial step-sizes $\alpha _0^i$, on ADBB is explored. Let ${\cal N}\left( {\theta,\;\; \xi } \right)$ to denote a normal distribution
with the mean vector $\theta$ and the covariance matrix $\xi$.
The total samples are $N=1000$, and we set $m=10$, $n=100$.
For each agent $i$, a half of sample vectors ${c_{ij}}$ are generated by independent and identically distributed $\mathcal{N}\left( {{{[2, - 2]}^{\top }},2I_n} \right)$ with label ${b_{ij}} =  + 1$, while the others are sample vectors ${c_{ij}}$ generated by independent and identically distributed $\mathcal{N}\left( {{{[-2, 2]}^{\top }},2I_n} \right)$ with label ${b_{ij}} =  - 1$.
\begin{figure}[htbp]
\begin{minipage}[h!]{0.5\linewidth}
{\includegraphics[height=5.5cm,width=6.5cm]{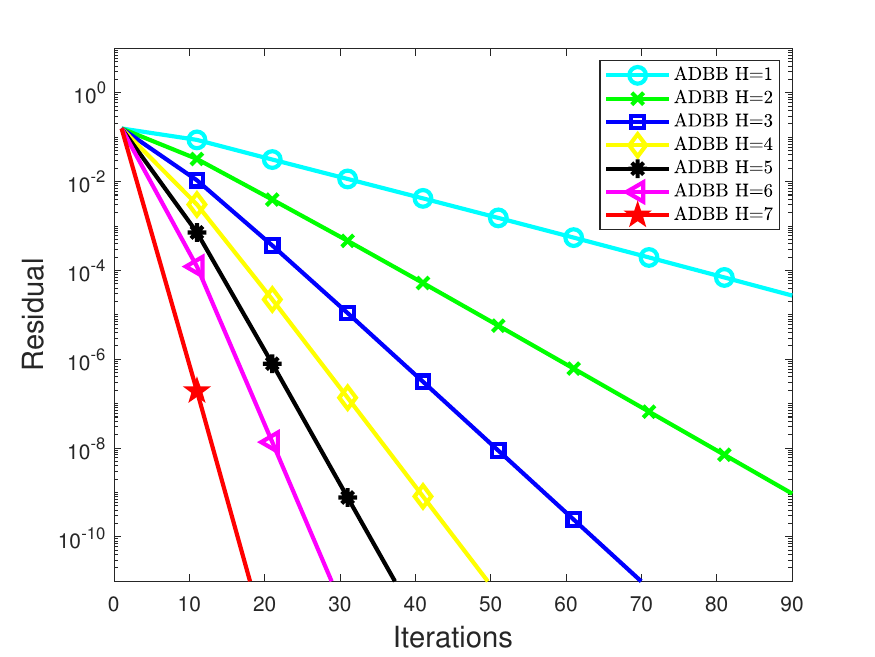}}
\caption{Performance of ADBB under different multi-consensus inner loop numbers $H$.}
\label{Fig.3}
\end{minipage}
\begin{minipage}[h!]{0.5\linewidth}
{\includegraphics[height=5.5cm,width=6.5cm]{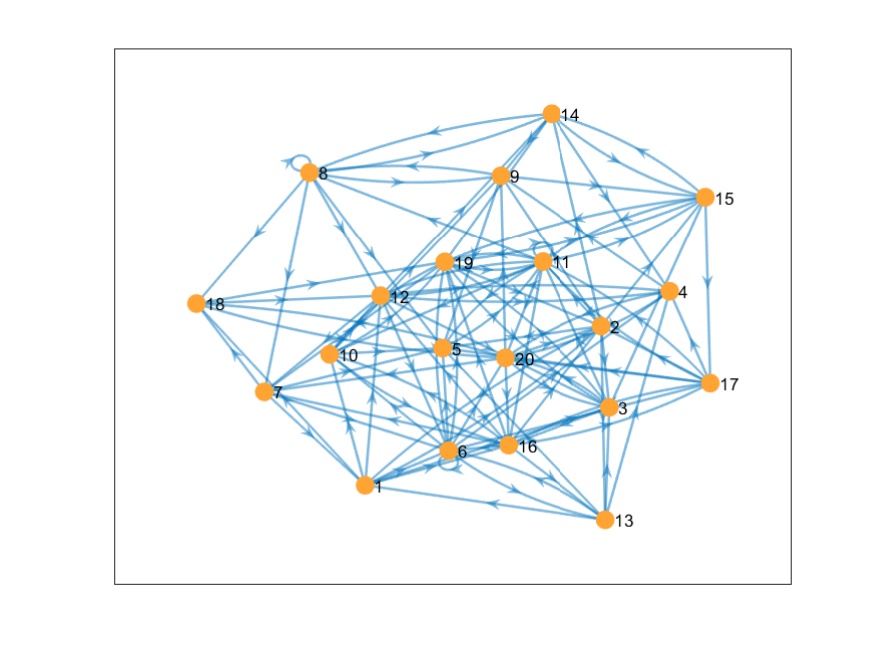}}
\caption{An unbalanced directed network with 20 agents.}
\label{Fig.4}
\end{minipage}
\end{figure}
\begin{figure}[h!]
\begin{minipage}[h!]{0.5\linewidth}
{\includegraphics[height=5.5cm,width=6.5cm]{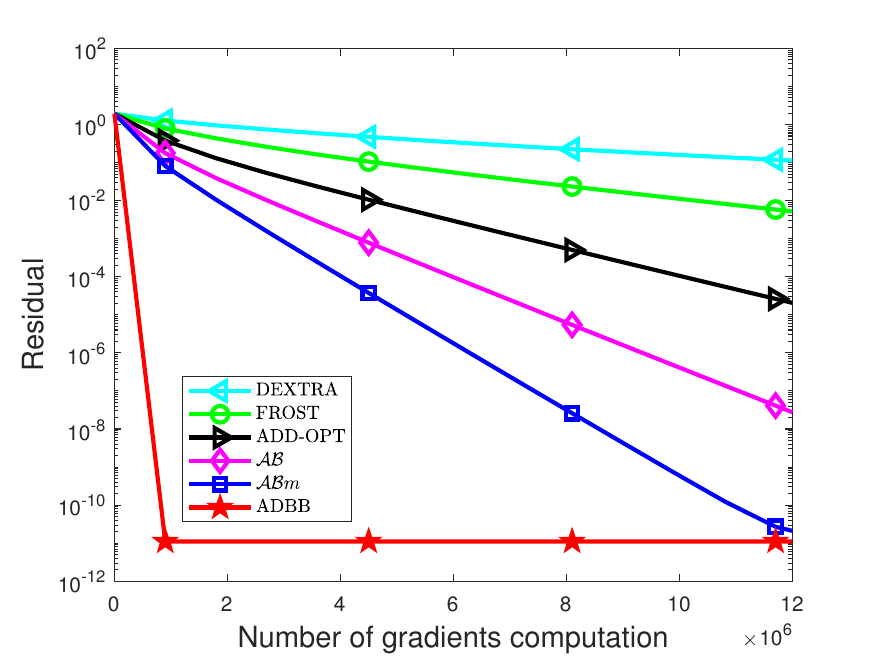}}
\caption{Performance comparison over gradients computation.}
\label{Fig.5}
\end{minipage}
\begin{minipage}[h!]{0.5\linewidth}
{\includegraphics[height=5.5cm,width=6.5cm]{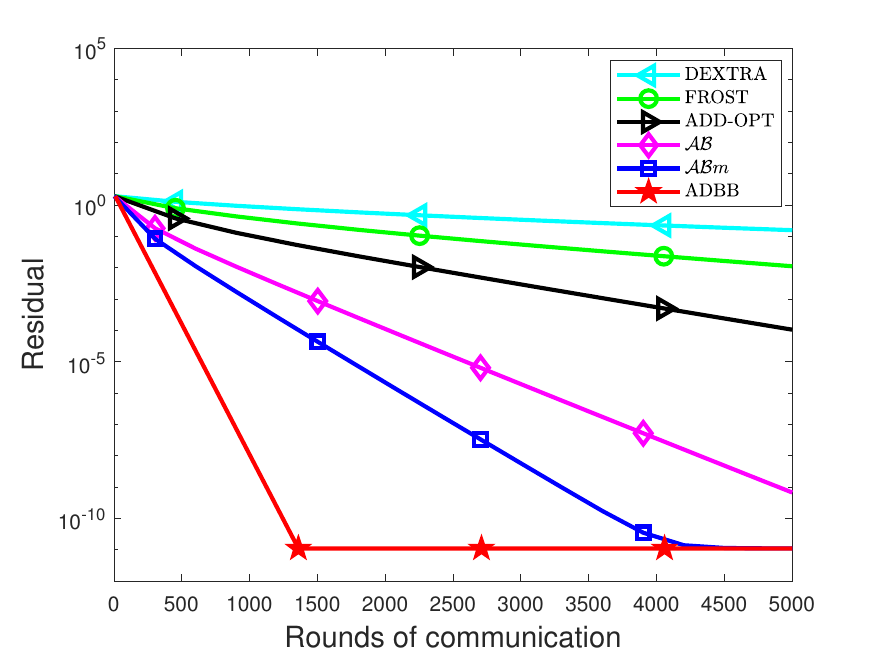}}
\caption{Performance comparison over rounds of communication.}
\label{Fig.6}
\end{minipage}
\end{figure}
Fig. \ref{Fig.1} shows the unbalanced directed communication network; Fig. \ref{Fig.2} compares the performance of ADBB with different initial step-sizes $\alpha _{\text{0}}^i$ by plotting the residual, where the multi-consensus inner loop number $H = 5$; Fig. \ref{Fig.3} compares the performance of ADBB with different multi-consensus inner loop number $H$ by plotting the residual, where the initial step-size $\alpha _{\text{0}}^i=1.2$. The results in Figs. \ref{Fig.2}-\ref{Fig.3} show that ADBB is not sensitive to initial step-sizes $\alpha _0^i$ and ADBB converges faster when the multi-consensus inner loop number $H$ increases, respectively.
\subsection*{Case Study 2}
In the second case, ADBB and some existing algorithms are used to identify whether a mushroom is poisonous or not according to its different features, such as cap-shape, cap-surface, cap-color, bruises, and so on. This case study is based on the mushroom data set provided in UCI Machine Learning Repository \cite{DuaD.andKarraTaniskidou2017}.
\begin{figure}[h!]
\begin{minipage}[h!]{0.5\linewidth}
{\includegraphics[height=5.5cm,width=6.5cm]{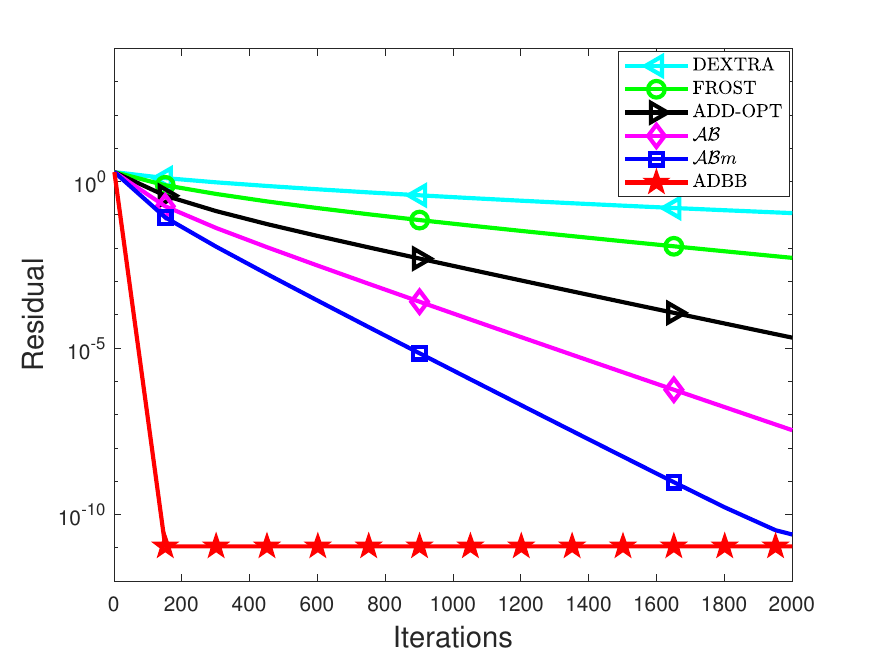}}
\caption{Performance comparison over iterations.}
\label{Fig.7}
\end{minipage}
\begin{minipage}[h!]{0.5\linewidth}
{\includegraphics[height=5.5cm,width=6.5cm]{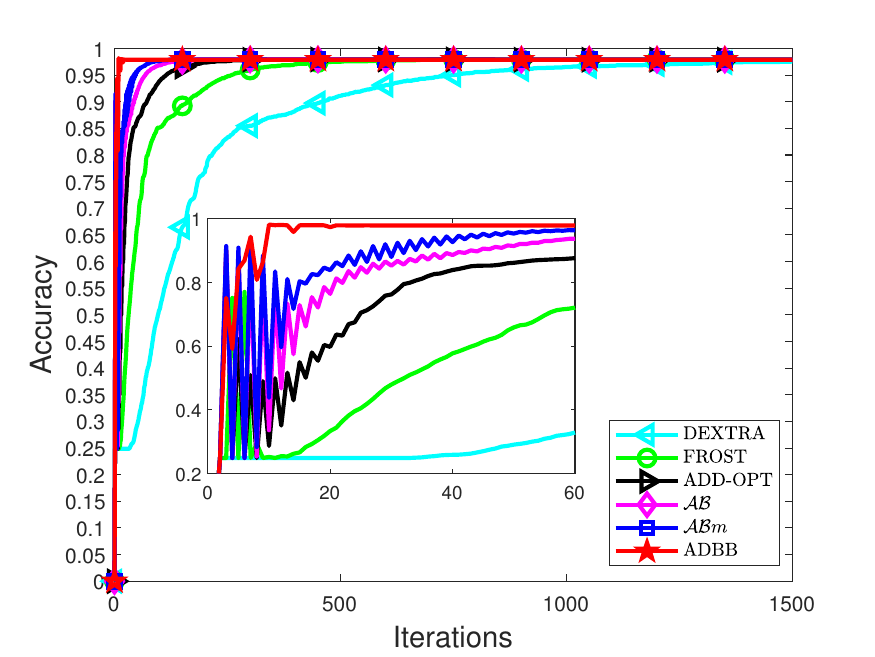}}
\caption{Testing accuracy rate.}
\label{Fig.8}
\end{minipage}
\end{figure}
We randomly choose 8124 samples from the data set, from which $N=6000$ samples are used to train the discriminator and the rest of samples are used for testing. Each sample has $n=112$ dimensions, which represents different features of the sample. We conduct the simulation in an unbalanced directed communication network as depicted in Fig. \ref{Fig.4}, where $m=20$; initial step-size $\alpha _0^i = 1.2$, $i \in \mathcal{V}$; multi-consensus inner loop number $H=3$. Assume label ${b_{ij}} =  + 1$ when the sample ${c_{ij}}$ is poisonous; Assume label ${b_{ij}} =  - 1$ when the sample ${c_{ij}}$ is edible.
Figs. \ref{Fig.5}-\ref{Fig.7} demonstrate that ADBB is competitive than most existing work over unbalanced directed networks in terms of gradients computation, rounds of communication, and the number of iterations.
The testing accuracy of different algorithms is plotted in Fig. \ref{Fig.8}. Tables \ref{Table 2}-\ref{Table 3} provide the confusion matrices to further clarify the testing results at iteration $k=100$ and $k=1000$, respectively. Notice that the total number of the testing samples is 2124 consisting of 1596 poisonous samples and 528 edible samples.

\begin{table}[!htbp]
\resizebox{120mm}{10mm}{
\begin{tabular}{|c|c|c|c|c|c|c|c|c|c|c|c|c|c|}
\hline
\multicolumn{2}{|c|}{ \multirow{2}*{\diagbox{\bf{True values}}{\bf{Predictive values}}} }& \multicolumn{2}{c|}{ADBB} & \multicolumn{2}{c|}{$\mathcal{A}\mathcal{B}m$}& \multicolumn{2}{c|}{$\mathcal{A}\mathcal{B}$}
 & \multicolumn{2}{c|}{ADD-OPT}& \multicolumn{2}{c|}{FROST}& \multicolumn{2}{c|}{DEXTRA}\\
\cline{3-14}
\multicolumn{2}{|c|}{}&P.&E.&P.&E.&P.&E.&P.&E.&P.&E.&P.&E.\\
\hline
\multicolumn{2}{|c|}{P.}&1587&9&1579&17&1563&33&1501&95&1297&299&553&1043\\
\cline{1-14}
\multicolumn{2}{|c|}{E.}&36&492&34&494&34&494&36&492&12&516&0&528\\
\cline{1-14}
\multicolumn{14}{|c|}{P. is the abbreviation of poisonous and E. is the abbreviation of edible}\\
\cline{1-14}
\end{tabular}}
\centering
\caption{The confusion matrix (Testing results at iteration $k=100$).}
\label{Table 2}
\end{table}
\begin{table}[!htbp]
\resizebox{120mm}{17mm}{
\begin{tabular}{|c|c|c|c|c|c|c|c|c|c|c|c|c|c|}
\hline
\multicolumn{2}{|c|}{ \multirow{2}*{\diagbox{\bf{True values}}{\bf{Predictive values}}} }& \multicolumn{2}{c|}{ADBB} & \multicolumn{2}{c|}{$\mathcal{A}\mathcal{B}m$}& \multicolumn{2}{c|}{$\mathcal{A}\mathcal{B}$}
 & \multicolumn{2}{c|}{ADD-OPT}& \multicolumn{2}{c|}{FROST}& \multicolumn{2}{c|}{DEXTRA}\\
\cline{3-14}
\multicolumn{2}{|c|}{}&P.&E.&P.&E.&P.&E.&P.&E.&P.&E.&P.&E.\\
\hline
\multicolumn{2}{|c|}{P.}&1589&7&1585&11&1579&17&1570&26&1564&32&1554&42\\
\cline{1-14}
\multicolumn{2}{|c|}{E.}&28&500&30&498&30&498&32&496&32&496&30&498\\
\cline{1-14}
\multicolumn{2}{|c|}{The best accuracy among 20 experiments}&\multicolumn{2}{c|}{0.98352}&\multicolumn{2}{c|}{0.98069}&\multicolumn{2}{c|}{0.97787}&\multicolumn{2}{c|}{0.97928}&\multicolumn{2}{c|}{0.97269}&\multicolumn{2}{c|}{0.96610}\\
\cline{1-14}
\multicolumn{2}{|c|}{The average accuracy among 20 experiments}&\multicolumn{2}{c|}{0.98021}&\multicolumn{2}{c|}{0.97712}&\multicolumn{2}{c|}{0.97710}&\multicolumn{2}{c|}{0.97628}&\multicolumn{2}{c|}{0.97124}&\multicolumn{2}{c|}{0.96573}\\
\cline{1-14}
\multicolumn{2}{|c|}{The worse accuracy among 20 experiments}&\multicolumn{2}{c|}{0.97012}&\multicolumn{2}{c|}{0.96922}&\multicolumn{2}{c|}{0.96975}&\multicolumn{2}{c|}{0.96828}&\multicolumn{2}{c|}{0.96823}&\multicolumn{2}{c|}{0.96521}\\
\cline{1-14}
\multicolumn{14}{|c|}{P. is the abbreviation of poisonous and E. is the abbreviation of edible}\\
\cline{1-14}
\end{tabular}}
\centering
\caption{The confusion matrix (Testing results at iteration $k=1000$).}
\label{Table 3}
\end{table}
\section{Conclusions and future work}\label{section six}
In this paper, a novel accelerated distributed algorithm constructing only row-stochastic weight matrices and using BB step-sizes, termed as ADBB, is developed to solve distributed convex optimization problems over an unbalanced directed network. Owing to the use of BB step-sizes and multi-consensus inner loops, ADBB allows each agent to automatically compute their step-sizes according to its local information and ensures the selection of larger step-sizes when achieving the globally optimal solution. Besides, ADBB has the fewer computation and communication costs than most existing distributed algorithms over unbalanced directed networks in simulations. To our knowledge, some existing variance-reduced stochastic gradient methods can protect the deterministic gradient algorithms from evaluating the full local gradients at each iteration and thus further reduce the computational costs. Therefore, as future work, it would be of considerable interest to incorporate the variance-reduced method into the proposed algorithm for distributed stochastic optimization over unbalanced directed networks.
\section*{Acknowledgments}
The work described in this paper was supported in part by the Fundamental Research Funds for the Central Universities under Grant XDJK2019AC001, in part by the Innovation Support Program for Chongqing Overseas Returnees under Grant cx2019005, and in part by the National Natural Science Foundation of China under Grant 61773321, Grant 61673080.
%
%
\bibliography{ADBB}
\end{document}